\documentclass{amsart}
\usepackage{amsmath,amsfonts,amsthm,enumerate}
\newcommand{\End}{\text{\rm End}}

\newcommand{\lemref}[1]{Lemma~\ref{#1}}
\newcommand{\eqnref}[1]{~(\ref{#1})}

\newtheorem{thm}{Theorem}[section]
\newtheorem{lem}[thm]{Lemma}

\theoremstyle{definition}

\numberwithin{equation}{section}

\title{}
\begin{document}
%\date{}

\title{Intermediate Wakimoto modules for Affine
$\mathfrak{sl}(n+1)$}
\author{Ben L.  Cox}
\author{Vyacheslav  Futorny}
%\runningauthor{Ben L. Cox and Vyacheslav Futorny}
%\runningtitle{Fock Space Realizations of Twisted Wakimoto Modules}
\address{Department of Mathematics \\
College of Charleston \\
66 George Street  \\
Charleston SC 29424, USA}\email{coxbl@cofc.edu}
%\date{April 18, 2001}

\address{School of Mathematics and Statistics\\
  University of Sydney\\
Sydney 2006, Australia\\
  On leave from Institute of Mathematics\\
University of S\~ao Paulo \\
Caixa Postal 66281  \\
S\~ao Paulo, CEP 05315-970, Brazil}\email{futorny@ime.usp.br}

\begin{abstract}  We construct certain boson type
  realizations of
affine $\mathfrak{sl}(n+1)$ that depend on a parameter $0\leq
r\leq n$ such that when $r=0$ we get a Fock space realization appearing
in the work of the first author and when $r=n$ they are the
Wakimoto modules described in the
work of Feigin and Frenkel.
   \end{abstract}
%\keywords{Wakimoto Modules, Intermediate Wakimoto Modules, Affine Lie
%Algebras, Fock Spaces}

%\classification{MR Primary}{17B67, 81R10}

\maketitle
\newpage
\section{Introduction}

Wakimoto modules for affine Lie algebras were introduced by B.~Feigin and
E.~Frenkel in 1988 by a homological characterization$^5$. These
modules  admit a remarkable boson realization on the Fock space
due to
Wakimoto$^{14}$ for $\hat{\mathfrak{sl}}(2)$,
and B.~Feigin and E.~Frenkel$^6$ for
$\hat{\mathfrak{sl}}(n)$ which plays an important role in the conformal
field theory providing a new bosonization rule for the Wess-Zumino-Witten
models.  Wakimoto modules have a geometric interpretation as certain
sheaves on a semi-infinite flag manifold described in  B.~Feigin and
E.~Frenkel$^{6}$. They  belong to the category $\mathcal O$  and
generically are isomorphic to corresponding Verma modules.  There are
numerous other authors who have explicitly constructed Wakimoto modules
for affine Lie algebras other than $\mathfrak{sl}(n+1)$.

Affine Lie algebras admit Verma type modules associated with
  non-standard Borel subalgebras which is described in the work of B.~
Cox$^3$, S.~Futorny and H.~Saifi$^9$ and H.~Jakobsen and V.~Kac$^{11}$.
In particular modules associated with the  {\it natural
Borel subalgebra} were first introduced by H.~Jakobsen and V.~Kac in
1985$^{11}$. They were studied by V.~
Futorny$^8$ under the name of  {\it imaginary Verma modules}.

A Fock space realization of the imaginary Verma modules for
$\hat{\mathfrak{sl}}(2)$ were constructed by Bernard and Felder$^1$
and then extended by the first author to the case of
$\hat{\mathfrak{sl}}(n)^4$. These realizations are given generically by
certain  Wakimoto type modules.

The main motivation for our work was a problem of finding  suitable
boson type realizations for all Verma type modules over 
$\hat{\mathfrak{sl}}(n+1)$.
In Theorem \ref{realization} we construct
   such realizations, {\it intermediate Wakimoto modules},
  for a series of generic Verma type modules depending on
  the parameter $0\leq r\leq n$. If $r=n$ this construction
coincides with the boson realization of Wakimoto modules in B.~Feigin and
E.~Frenkel$^5$.   On the other hand when $r=0$ the obtained
representation gives a Fock space  realization described in the work of
the first author$^4$. Using this realization we plan to discuss the
detailed structure of intermediate Wakimoto modules in a subsequent
paper.

\section{Preliminaries}

Fix a positive  integer
$n$, $0\leq r\leq n$,
$\gamma\in\mathbb C^*$. Set
$k=\gamma^2-(r+1)$.
Let $\mathfrak g=\mathfrak{sl}(n+1,\mathbb C)$ and let
  $E_{ij}$, $i,j=1, \ldots, n+1$ be
the standard basis for  $\mathfrak{gl}(n+1,\mathbb C)$.  Set
$H_i:=E_{ii}-E_{i+1,i+1}$, $E_i:=E_{i,i+1}$, $F_i:=E_{i+1,i}$ which is a
basis for $\mathfrak{sl}(n+1,\mathbb C)$. Furthermore
we denote the Killing form by
$(X|Y)=\text{tr}\,(XY)$ and $X_m=t^m\otimes X$ for
$X,Y\in \mathfrak g$ and $m\in\mathbb Z$.  Let
$\{\alpha_1,\dots,\alpha_n\}$ be a base for $\Delta^+$, the positive set
of
roots for $\mathfrak g$, such that $H_i=\check\alpha_i$
and let
$\Delta_r$ be the root system with base $\{\alpha_1,\dots,\alpha_r\}$
($\Delta_r=\emptyset$, if $r=0$) of the Lie subalgebra $\mathfrak g_r=
\mathfrak{sl}(r+1,\mathbb C)$. A Cartan subalgebra
$\mathfrak H$ (respectively $\mathfrak H_r$) of $\mathfrak g$
(respectively $\mathfrak g_r$) is spanned by $H_i$, $i=1, \ldots, n$
(respectively $i=1, \ldots, r$) and set $\mathfrak H_0=0$.

  For any Lie algebra $\mathfrak a$, let
$L(\mathfrak a
)=\mathbb{C} [t,t^{-1}] \bigotimes\mathfrak a$ be the loop algebra of
$\mathfrak a$. Then   $\hat{\mathfrak g}=\hat{\mathfrak{sl}}(n+1,\mathbb 
C) =L({\mathfrak g}
)\oplus\mathbb{C}
c\oplus\mathbb{C} d$ and $\hat{\mathfrak g}_r=L({\mathfrak g_r}
)\oplus\mathbb{C}
c\oplus\mathbb{C} d$
  are the associated  affine Kac-Moody
algebras with $\hat{\mathfrak H}=\mathfrak H\oplus\mathbb{C}
c\oplus\mathbb{C} d$ and $\hat{\mathfrak H}_r=\mathfrak 
H_r\oplus\mathbb{C}
c\oplus\mathbb{C} d$ respectively.

The  algebra
$\hat{\mathfrak g}$ has generators
$E_{im}, F_{im}, H_{im}$,  $i=1, \ldots, n$,  $m\in\mathbb Z$, and 
central element $c$ with the product
$$[X_m, Y_n]=t^{m+n}[X,Y]+ \delta_{m+n,0}m(X|Y)c.$$

%\section{Notation and Realization}

\subsection{ Oscillator algebras}

Let $\hat{\mathfrak a}$ be the infinite
dimensional Heisenberg algebra with generators
$a_{ij,m}$, $a_{ij,m}^*$, and $\mathbf 1$,  $1\leq i\leq j\leq n$  and
$m\in
\mathbb Z$, subject to the relations
\begin{align*}
[a_{ij,m},a_{kl,n}]&=[a_{ij,m}^*,a^*_{kl,n}]=0, \\
[a_{ij,m},a^*_{kl,n}]&=\delta_{ik}\delta_{jl}
\delta_{m+n,0}\mathbf 1, \\
[a_{ij,m},\mathbf 1]&=[a^*_{ij,m},\mathbf 1]=0.
\end{align*}
Such an algebra has a representation $\tilde{\rho}:\hat{\mathfrak
a}\to {\mathfrak{gl}}(\mathbb
C[\mathbf x])$ where
\begin{align*}   \mathbb C[\mathbf x]&:=
        \mathbb C[x_{ij,m}|i,j,m\in \mathbb Z,\,1\leq i\leq j\leq
     n]
\end{align*}
denotes the algebra over $\mathbb C$ generated by the indeterminates
$x_{ij,m}$ and $\tilde{\rho}$ is
defined by
\begin{align*}
  \tilde{\rho}( a_{ij,m}):&=\begin{cases}
   \partial/\partial
x_{ij,m}&\quad \text{if}\quad m\geq 0,\enspace\text{and}\enspace  j\leq r
\\ x_{ij,m} &\quad \text{otherwise},
\end{cases}
  \\
\tilde{\rho}(a_{ij,m}^*):&=
\begin{cases}x_{ij,-m} &\enspace \text{if}\quad m\leq
0,\enspace\text{and}\enspace j\leq r \\ -\partial/\partial
x_{ij,-m}&\enspace \text{otherwise}. \end{cases}
\end{align*}
and $\tilde{\rho}(\mathbf 1)=1$.
In this case
$\mathbb C[\mathbf x]$ is an
$\hat{\mathfrak a}$-module generated by $1=:|0\rangle$, where
$$
a_{ij,m}|0\rangle=0,\quad m\geq  0 \enspace\text{and}\enspace j\leq r,
\quad a_{ij,m}^*|0\rangle=0,\quad m>0\enspace\text{or}\enspace j>r.
$$
Let $\hat{\mathfrak a}_r$ denote the subalgebra generated by $a_{ij,m}$
and
$a_{ij,m}^*$ and $\mathbf 1$, where $1\leq i\leq j\leq r$ and
$m\in\mathbb
Z$. If $r=0$, we set $\hat{\mathfrak a}_r=0$.

Let $A_n=((\alpha_i|\alpha_j))$ be the Cartan matrix for 
$\mathfrak{sl}(n+1,\mathbb C)$
and
let $\mathfrak B$ be the matrix whose entries are
$$
\mathfrak B_{ij}:=(\alpha_i|\alpha_j)(\gamma^2
     -\delta_{i>r}\delta_{j>r}(r+1)
     +\frac{r}{2}\delta_{i,r+1}\delta_{j,r+1})
$$

where $$
\delta_{i>r}= \begin{cases}
1& \quad \text{if}\quad i>r,\\
0& \quad \text{otherwise}.
\end{cases}$$

In other words
$$
\mathfrak B:=\gamma^2A_n -(r+1)
\begin{pmatrix} 0  & 0 \\ 0 & A_{n-r} \end{pmatrix}+rE_{r+1,r+1}.
$$

We also have the Heisenberg Lie algebra $\hat{\mathfrak b}$ with 
generators $b_{im}$,
$1\leq
i\leq n$, $m\in\mathbb Z$, $\mathbf 1$, and relations
$[b_{im},b_{jp}]=m\,\mathfrak B_{ij}\delta_{m+p,0}\mathbf 1$ and
$[b_{im},\mathbf 1]=0$.

  For each
$1\leq i\leq n$ fix $\lambda_i\in\mathbb C$ and let
$\lambda=(\lambda_1, \ldots, \lambda_n)$.
  Then the algebra  $\hat{\mathfrak b}$
has
a representation
$\rho_{\lambda}:\hat{\mathfrak b}\to\End(\mathbb C[\mathbf y]_\lambda)$ 
where
\begin{align*}   \mathbb C[\mathbf y]&:=
        \mathbb C[y_{i,m}|i,m\in \mathbb N^*,\,1\leq i\leq
     n]
\end{align*}
   and $\rho_{\lambda}$ is
defined
on
$\mathbb C[\mathbf y]$ defined by
$$
\rho_{\lambda}(b_{i0})= \lambda_i, \quad \rho_{\lambda}(b_{i,-m})=
      \mathbf e_i\cdot \mathbf y_{m},\quad
      \rho_{\lambda}(b_{im})=m\mathbf e_i \cdot  
\,\frac{\partial}{\partial
\mathbf y_{m}}\quad
\text{for}\quad m>0
$$
and $\rho_{\lambda}(\mathbf 1)=1$.  Here
$$
\mathbf y_m=(y_{1m},\cdots, y_{nm}),\quad \frac{\partial}{\partial
\mathbf y_{m}}=\left(\frac{\partial}{\partial   y_{1m}}, \cdots,
\frac{\partial}{\partial y_{nm}} \right)
$$
and $\mathbf e_i$ are vectors in $\mathbb C^n$ such that $\mathbf
e_i\cdot
\mathbf e_j=\mathfrak B_{ij}$ where $\cdot $ means the usual dot
product.

Note that since $\mathfrak B_{ij}$ is symmetric, it is
orthogonally diagonalizable, (i.e. there exists an orthogonal matrix $P$
such that $P^t\mathfrak BP$ is a diagonal matrix) and hence we can find
vectors $\mathbf e_i$  in
$\mathbb C^n$ such that
$\mathbf e_i\cdot
\mathbf e_j=\mathfrak B_{ij}$. In fact for
$m>0$ and
$n<0$ we get
\begin{align*}
[b_{im},b_{jn}]&=[m\mathbf e_i \cdot  \,\frac{\partial}{\partial
\mathbf y_{m}},\mathbf e_j\cdot \mathbf y_{-n}]  \\
&=m\sum_{k,l}[  e_{ik} \frac{\partial}{\partial
y_{km}},  e_{jl}  y_{l,-n}]  \\
&=m\delta_{m+n,0}\sum_{k} e_{ik}  e_{jk}  =m\delta_{m+n,0}\mathfrak
B_{ij}.
\end{align*}
(See also  the work of B. Feigin and E. Frenkel$^7$.)

\subsection{Formal Distributions}
We need some more notation that will simplify some of the arguments
later.
This notation follows the books of A. Matsuo and K. Nagatomo$^{14}$ and 
V.
Kac$^{12}$:  A {\it formal distribution} is an expression of the form
$$
a(z,w,\dots)=\sum_{m,n,\dots\in\mathbb Z}a_{m,n,\dots}z^mw^n
$$
where the $a_{m,n,\dots}$ lie in some fixed vector space $V$.
We define $\partial a(z)=\partial_z a(z)=\sum_nna_nz^{n-1}$.
We also have expansion about zero:  there are two
canonical embeddings of fields $\iota_{z,w}:\mathbb C(z-w)\to
\mathbb C[[z,w]]$ and  $\iota_{w,z}:\mathbb C(z-w)\to
\mathbb C[[z,w]]$
where $
\iota_{z,w}(a(z,w))$ is formal Laurent series expansion in $z^{-1}$
and
$-\iota_{w,z}(a(z,w))$  is formal Laurent series expansion in $z$.
The {\it formal delta function}
$\delta(z-w)$ is the formal distribution
$$
\delta(z-w)=z^{-1}\sum_{n\in\mathbb Z}\left(\frac{z}{w}\right)^n
=\iota_{z,w}\left(\frac{1}{z-w}\right)-\iota_{w,z}\left(
\frac{1}{z-w}\right).
$$
For any sequence of elements $\{a_{(m)}\}_{m\in
\mathbb Z}$ in the ring $\End (V)$, $V$ a vector space,  the
formal distribution
\begin{align*}
a(z):&=\sum_{m\in\mathbb Z}a_{(m)}z^{-m-1}
\end{align*}
is called a {\it field}, if for any $v\in V$, $a_{(m)}v=0$ for
$m\gg0$. For a field such that $a_{(m)}$ are creation operators for $m\ll
0$, we set
\begin{align*}
     a^-(z):&=\sum_{m\geq 0}a_{(m)}z^{-m-1},\quad\text{and}\quad
    a^+(z):=\sum_{m<0}a_{(m)}z^{-m-1}.
\end{align*}
Observe that $a_{ij}(z)$ for $j>r$ is not a field whereas $a_{ij}^*(z)$
is always a field.
We also define
$$
\delta^-(z-w)=\iota_{z,w}\left(\frac{1}{z-w}\right),\quad
\delta^+(z-w)=-\iota_{w,z}\left(\frac{1}{z-w}\right).
$$
Note that
$$
-\partial_z\delta(z-w)=\partial_w\delta(z-
w)=\iota_{z,w}\left(\frac{1}{(z-w)^2}\right)
-\iota_{w,z}\left(\frac{1}{(z-w)^2}\right).
$$
Finally we use the convention
\begin{align*}
     a^-_{ij}(z):&=0,\quad\text{and}\quad
    a^+_{ij}(z):=a_{ij}(z),
\end{align*}
\begin{align*}
     a^{*-}_{ij}(z):&=a^{*}_{ij}(z),\quad\text{and}\quad
    a^{*+}_{ij}(z):=0\quad \text{for}\quad j>r.
\end{align*}

  The {\it normal ordered product} of two formal distributions
$a(z)$ and
$b(w)$ is
defined by
$$
:a(z)b(w):=a^+(z)b(w)+b(w)a^-(z).
$$

For any $1\leq i\leq j\leq n$, we define
$$
a^*_{ij}(z)=\sum_{n\in\mathbb Z}a^*_{ij,n}z^{-n},\quad
a_{ij}(z)=\sum_{n\in\mathbb Z}a_{ij,n}z^{-n-1}
$$
and
$$
b_i(z)=\sum_{n\in\mathbb Z}b_{in}z^{-n-1}.
$$

In this case

\begin{align*}
[b_i(z),b_j(w)]
%=\sum_{m,n}[b_{i,m},b_{j,n}]z^{-m-1}w^{-n-1}
%=\sum_{m,n}m\delta_{i,j}\delta_{m+n,0}z^{-m-1}w^{-n-1}
%=\sum_{m}m\delta_{i,j}z^{-m-1}w^{m-1}
&=\mathfrak{B}_{ij}\partial_w\delta(z-w), \\
[a_{ij}(z),a^*_{kl}(w)]&
%   =\sum_{m,n}[a_{\alpha,m},a^*_{\beta,n}]z^{-m-1}w^{-n}
%   =\sum_{m}\delta_{\alpha,\beta}
%     \mathbf 1z^{-m-1}w^{m}
     =\delta_{ik}\delta_{jl}\mathbf 1\delta(z-w).
\end{align*}

Let
\begin{equation}
\lfloor
ab\rfloor=a(z)b(w)-:a(z)b(w):=[a^-(z),b(w]
\end{equation}
(half of
$[a(z),b(w)]$) denote the {\it contraction} of any two formal 
distributions
$a(z)$ and $b(w)$. For example if $j,l\leq r$, then
\begin{align}
\lfloor a_{ij} a_{kl}^*\rfloor
&%=\sum_{m\geq 0}\delta_{ik}\delta_{jl}\delta_{m+n,0}z^{-m-1}
%w^{-n}
=\sum_{m\geq 0}\delta_{ik}\delta_{jl}z^{-m-1}w^{m}
=\delta_{i,k}\delta_{j,l}\delta^-(z-w)
=\delta_{ik}\delta_{jl}
\,\iota_{z,w}\left(\frac{1}{z-w}\right)\\
\lfloor a_{kl}^*a_{ij}\rfloor
&
%=-\sum_{m\geq 1}\delta_{ik}\delta_{jl}\delta_{m+n,0}z^{-m}
%w^{-n-1}
=-\sum_{n<0}\delta_{ik}\delta_{jl}z^{n}w^{-n-1}
=-\delta_{ik}\delta_{jl}\delta^+(w-z)=
\delta_{ik}\delta_{jl}\,\iota_{z,w}\left(\frac{1}{w-z}
\right).
\end{align}

We need the very useful Wick's Theorem$^{2,12,14}$:
\begin{thm}  Let $a_i(z)$
and
$b_j(z)$ be formal distributions with coefficients in the associative
algebra
  $\End(\mathbb C[\mathbf x]\otimes \mathbb C[\mathbf y])$,
  satisfying
\begin{enumerate}
\item $[ \lfloor a_ib_j\rfloor ,c_k(z)]=0$, for all $i,j,k$ and $c=a$ or

$c=b$,
\item $[a^\pm_i(z),b^\pm_j(w)]=0$ for all $i$ and $j$.
\end{enumerate}
Then
$$
:a_1(z)a_2(z)\cdots a_k(z):b(w)=\sum_{i=1}^k:a_1(z)\cdots \lfloor
a_ib\rfloor\cdots a_k(z):
$$
and
\begin{align*}
:&a_1(z)\cdots a_m(z)::b_1(w)\cdots
b_k(w):= \\
   &\sum_{s=0}^{\min(m,k)}\sum_{\genfrac{}{}{0pt}{}{i_1<\cdots<i_s,}
{j_1\neq \cdots \neq j_s}}\lfloor a_{i_1}b_{j_1}\rfloor\cdots
\lfloor a_{i_s}b_{i_s}\rfloor
:a_1(z)\cdots a_m(z)b_1(w)\cdots
b_k(w):_{(i_1,\dots,i_s;j_1,\dots,j_s)}
\end{align*}
where the subscript ${(i_1,\dots,i_s;j_1,\dots,j_s)}$ means that
those factors $a_i(z)$, $b_j(w)$ with indices
$i\in\{i_1,\dots,i_s\}$, $j\in\{j_1,\dots,j_s\}$ are to be omitted from
the product \newline
$:a_1(z)\cdots a_m(z)b_1(w)\cdots b_k(w):$.
\end{thm}

The proof is identical to that in Kac$^{12}$ even though it is stated
for fields $a_i(z)$ and $b_j(z)$ in that text.

We will also need the following two
results.
\begin{thm}[Taylor's Theorem$^{12,14}$]
\label{Taylorsthm}  Let
$a(z)$ be a formal distribution.  Then in the region $|z-w|<|w|$,
\begin{equation}
a(z)=\sum_{j=0}^\infty \partial_w^{(j)}a(w)(z-w)^j.
\end{equation}
\end{thm}

\begin{thm}[Kac$^{12}$, Theorem 2.3.2]\label{ope} Let $a(z)$ and $b(z)$
be the formal distributions with coefficients in the associative algebra
  $\End(\mathbb C[\mathbf x]\otimes \mathbb C[\mathbf y])$.  The following
are equivalent
\begin{enumerate}[(i)]
\item
$\displaystyle{[a(z),b(w)]=\sum_{j=0}^{N-1}\partial_w^{(j)}
\delta(z-w)c^j(w)}$, where $c^j(w)$ is a formal distribution with
coefficients in
  $\End(\mathbb C[\mathbf x]\otimes \mathbb C[\mathbf y])$.
\item
$\displaystyle{\lfloor
ab\rfloor=\sum_{j=0}^{N-1}\iota_{z,w}\left(\frac{1}{(z-w)^{j+1}}
\right)
c^j(w)}$.
\end{enumerate}\label{Kacsthm}
\end{thm}

In other words the singular part of the {\it operator product
expansion}
$$
\lfloor
ab\rfloor=\sum_{j=0}^{N-1}\iota_{z,w}\left(\frac{1}{(z-w)^{j+1}}
\right)c^j(w)
$$
completely determines the bracket of mutually local formal
distributions $a(z)$ and $b(w)$.   One writes
$$
a(z)b(w)\sim \sum_{j=0}^{N-1}\frac{c^j(w)}{(z-w)^{j+1}}.
$$
For example
$$
b_i(z)b_j(w)\sim \frac{\delta_{ij}}{(z-w)^2}.
$$

\subsection{Verma type modules}

For a Lie algebra $\mathfrak a$ we denote by $U(\mathfrak a)$ the 
universal enveloping algebra of $\mathfrak a$.

   Let
  ${\mathfrak g}_{\alpha}$ be a root subspace of $\mathfrak g$ 
corresponding to a root
$\alpha$,
  ${\mathfrak
n}^{\pm}=\oplus_{\alpha\in {\Delta}^+}{\mathfrak g}_{\pm \alpha}$ and
  $\mathfrak g={\mathfrak{n}}^-\oplus\mathfrak H\oplus
{\mathfrak{n}}^+$ a Cartan decomposition of $\mathfrak g$. Denote also
  ${\mathfrak
n}^{\pm}_r={\mathfrak n}^{\pm}\cap {\mathfrak g}_r$,
${\mathfrak n}^+(r)={\mathfrak n}^+\setminus {\mathfrak n}^+_r$,
$$\bar{B}_r=L({\mathfrak n}^+(r))\oplus ({\mathfrak n}^+_r\otimes
\mathbb C[t]) \oplus (({\mathfrak n}^-_r)\oplus {\mathfrak H})\otimes
\mathbb C[t]t).$$
  Then
  $B_r=\bar{B}_r\oplus \hat{\mathfrak H}$ is a Borel subalgebra of
$\hat{\mathfrak g}$ for any $0\leq r\leq n$.

Fix $\tilde{\lambda}\in \hat{\mathfrak H}^*$ and consider a
$\hat{\mathfrak g}$-module

$$M_r(\tilde{\lambda})=U(\hat{\mathfrak g})\otimes_{U(B_r)}\mathbb C 
v_{\tilde{\lambda}}$$
where $\bar{B}_r  v_{\tilde{\lambda}}=0$ and $h v_{\tilde{\lambda}}=
\tilde{\lambda}(h) v_{\tilde{\lambda}}$ for all $h\in \hat{\mathfrak H}$.

Module  $M_r(\tilde{\lambda})$ is a particular case of a Verma type
module studied in Cox$^3$, Futorny and Saifi$^9$. When $r=n$ it
gives a usual Verma module  construction. If $r=0$ we get an imaginary
Verma module.

Let $\tilde{\lambda}_r=\tilde{\lambda}|_{\hat{\mathfrak H}_r}$.
Verma type module  $M_r(\tilde{\lambda})$ contains a
$\hat{\mathfrak g}_r$-submodule $M(\tilde{\lambda}_r)=U(\hat{\mathfrak 
g}_r)(1\otimes  v_{\tilde{\lambda}})$ which is isomorphic to a usual 
Verma module for
$\hat{\mathfrak g}_r$.

Note that the proof given in Kac's book cited above works also in the
setting that the distributions are not necessarily fields.

\begin{thm}[Cox$^3$, Futorny and Saifi$^9$]\label{vermatype}
Let $\tilde{\lambda}(c)\neq 0$. Then the submodule structure of
$M_r(\tilde{\lambda})$ is completely determined by the submodule 
structure of
$M(\tilde{\lambda}_r)$. In particular, $M_r(\tilde{\lambda})$ is 
irreducible if
$M(\tilde{\lambda}_r)$ is irreducible.

\end{thm}

%\bigskip

\section{Intermediate Wakimoto modules}

Define for $\quad 1\leq i\leq n$,
$$
E_{i}(z)=\sum_{n\in\mathbb Z}E_{in}z^{-n-1},\quad
F_{i}(z)=\sum_{n\in\mathbb Z}F_{in}z^{-n-1},\quad
H_i(z)=\sum_{n\in\mathbb
Z} H_{in}z^{-n-1}.
$$
The defining relations between the generators
of $\hat{\mathfrak g}$ can be written as follows

\begin{align*}
[H_i(z),H_j(w)]&=(\alpha_i|\alpha_j)c\partial_w \delta(w-z)
\tag{R1}\\
    [H_i(z),E_j(w)]&=(\alpha_i|\alpha_j)E_j(z)\delta(w-z)\tag{R2} \\
    [H_i(z),F_j(w)]&=-(\alpha_i|\alpha_j)F_j(z)\delta(w-z)\tag{R3} \\
    [E_i(z),F_j(w)]&=\delta_{i,j}(H_i(z)\delta(w-z)+c\partial_w
      \delta(w-z))\tag{R4} \\
    [F_i(z),F_j(w)]&=[E_i(z),E_j(w)]=0\quad\text{if}\quad
(\alpha_i|\alpha_j)\neq -1\tag{R5} \\
    [F_i(z_1),F_i(z_2),F_j(w)]&=[E_i(z_1),E_i(z_2),E_j(w)]=0
\quad\text{if}\quad (\alpha_i|\alpha_j)=
    -1\tag{R6}
    \end{align*}
where $[X,Y,Z]:=[X,[Y,Z]]$ is the Engel bracket for any three
operators
$X,Y,Z$.

Recall that $\mathbb C[\mathbf x]$ is an $\hat{\mathfrak a}$-module with 
respect to
the representation $\tilde{\rho}$ and   $\mathbb C[\mathbf y]$ is a
$\hat{\mathfrak b}$-module with respect to
  ${\rho}_{\lambda}$. The main result of the paper is the following 
theorem where we define a representation
$$\rho:\hat{\mathfrak g}\rightarrow
\mathfrak{gl}(\mathbb C[\mathbf x]\otimes \mathbb C[\mathbf y]). $$
We  use the notation
$\rho(X_m):=\rho(X)_m$, for
$X\in\mathfrak g$.
%\vfill\eject

\begin{thm}\label{realization}
Let $\lambda\in \mathfrak H^*$ and set $\lambda_i=\lambda(H_i)$.  The
generating functions
\begin{align*}
\rho(c)&=\gamma^2-(r+1), \\
\rho(F_{i})(z)&=a_{ii}+\sum_{j=i+1}^{n}a_{ij}a_{i+1,j}^*,\\ \\
\rho(H_i)(z)&=2:a_{ii}a_{ii}^*:+\sum_{j=1}^{i-1}\left(:
          a_{ji}a_{ji}^*:-:a_{j,i-1} a_{j,i-1}^*:\right) \\
        &\quad+\sum_{j=i+1}^{n}\left(:a_{ij}a_{ij}^*:
         - :a_{i+1,j}a_{i+1,j}^*:\right)
          + b_i,\\ \\
\rho(E_i)(z)&= :a_{ii}^*\left(
        \sum_{k=1}^{i-1} a_{k,i-1}a_{k,i-1}^*-\sum_{k=1}^{i}
        a_{ki}a^*_{ki}\right): +\sum_{k=i+1}^n a_{i+1,k}a_{ik}^*
          -\sum_{k=1}^{i-1}a_{k,i-1}a_{ki}^* \\
&\quad-a_{ii}^*b_{i}-
      \left(\delta_{i> r}(r+1)+\delta_{i\leq
r}(i+1)-\gamma^2\right)\partial a_{ii}^* ,
\end{align*}
define an action of the generators
$E_{im}$, $F_{im}$,
$H_{im}$, $i=1, \ldots, n$, $m\in \mathbb Z$ and $c$, on the Fock space
$\mathbb C[\mathbf x]\otimes
\mathbb C[\mathbf y]$.  In the above $a_{ij}$, $a_{ij}^*$ and $b_i$ 
denotes
$a_{ij}(z)$,
$a_{ij}^*(z)$ and $b_i(z)$ respectively.
\end{thm}

Theorem \ref{realization} defines a boson type realization of
$\hat{\mathfrak{sl}}(n+1)$ and a module structure on the Fock
space
$\mathbb C[\mathbf x]\otimes \mathbb C[\mathbf y]$ that depends on the
parameter $r$, $0\leq r\leq n$.
We will call such a module, an {\it intermediate Wakimoto module} and
denote it by
$W_{n,r}(\lambda, \gamma)$.  The intermediate Wakimoto
modules $W_{n,r}(\lambda, \gamma)$ have the property that the subalgebra
$\bar{B}_{r}$ annihilates the vector $1\otimes 1\in \mathbb C[\mathbf 
x]\otimes \mathbb C[\mathbf y]$, $h(1\otimes 1)=\lambda(h)(1\otimes 1)$ 
for all $h\in \mathfrak H$ and
$c(1\otimes 1)=(\gamma^2 - (r+1))(1\otimes 1)$.
Consider the $\hat{\mathfrak g}_r$-submodule $W=U(\hat{\mathfrak
g}_r)(1\otimes 1)
\simeq W_{r,r}(\lambda, \gamma)$ of $W_{n,r}(\lambda, \gamma)$. Then
$W$ is isomorphic to the Wakimoto module $W_{\lambda(r),\tilde{\gamma}}$
  of Feigin and Frenkel$^{5}$ where $\lambda(r)=\lambda|_{\mathfrak
H_r}$,
$\tilde{\gamma}=\gamma^2-(r+1)$.

Consider $\tilde{\lambda}\in \hat{\mathfrak H}^*$ such that 
$\tilde{\lambda}|_{\mathfrak H}=\lambda$, 
$\tilde{\lambda}(c)=\gamma^2-(r+1)$, a Verma type module
$M_r(\tilde{\lambda})$ and its $\hat{\mathfrak g}_r$-submodule
$M(\tilde{\lambda}_r)$. Suppose that $M(\tilde{\lambda}_r)$ is 
irreducible. In this case the Wakimoto module  
$W_{\lambda(r),\tilde{\gamma}}$ is isomorphic
to $M(\tilde{\lambda}_r)$. Let
$\tilde W=U(\hat{\mathfrak g})W_{\lambda(r),\tilde{\gamma}}$ and assume 
that
$\lambda(c)\neq 0$. Then by Theorem \ref{vermatype} the module
$M_r(\tilde{\lambda})$ is irreducible and therefore isomorphic to
$\tilde W$. Hence Theorem \ref{realization} provides a boson type
  realization for generic Verma type modules.

We believe that generically Verma type modules and intermediate Wakimoto 
modules are isomorphic. A similar realization must exist for all Verma 
type modules over
$\hat{\mathfrak{sl}}(n+1)$ and other affine Lie algebras.

\bigskip

\section{Formal distribution computations}

Set
{
\begin{align*}
\mathcal H_{i}(z):&=2:a_{ii}a_{ii}^*:+\sum_{j=1}^{i-1}\left(:
          a_{ji}a_{ji}^*:-:a_{j,i-1} a_{j,i-1}^*:\right) \\
        &\quad+\sum_{j=i+1}^{n}\left(:a_{ij}a_{ij}^*:
         - :a_{i+1,j}a_{i+1,j}^*:\right).
%&=2(-1)^{i}\sum_{k\in\mathbb
%Z,n\in\mathbb
%N}:a_{\alpha_0+n\delta}(m-k)
%a^*_{\alpha_0+n\delta}(k):-
%:a_{\alpha_1+n\delta}(m-k)
%a^*_{\alpha_1+n\delta}(k):
\end{align*}}
   In the above $a_{ij}$, and $a_{ij}^*$  denotes
$a_{ij}(z)$,
$a_{ij}^*(z)$ respectively.

For any $\alpha\in\Delta^+$ we can find unique $1\leq k\leq l\leq n$
such
that
\begin{equation}\label{alphakl}
\alpha_{kl}:=\alpha=\alpha_k+\cdots +\alpha_l.
\end{equation}

Set $a_{\alpha}:=a_{kl}$ and $a_{\alpha}^*:=a_{kl}^*$.
Observe that

\begin{align*}
(\alpha_i|\alpha)=\sum_{j=k}^l(\alpha_i|\alpha_j)&=\left(2\delta_{ik}\delta_{
il}+
\delta_{k<i}\left(\delta_{il}-
    \delta_{i-1,l}\right)+\delta_{l> i}\left(\delta_{ik}-
    \delta_{i+1,k}\right)\right) \\
&=\delta_{ik}-\delta_{i+1,k}-\delta_{i,l+1}+\delta_{i+1,l+1}.
\end{align*}

Since this is the case we can rewrite
\begin{align*}
\mathcal H_{i}(z):&=\sum_{\alpha\in
\Delta^+}(\alpha_i|\alpha):a_{\alpha}a_{\alpha}^*:.
\end{align*}

Moreover we have
\begin{align}
\lfloor a_{\alpha} a_{\beta}^*\rfloor
&%=\sum_{m\geq 0}\delta_{i,k}\delta_{j,l}\delta_{m+n,0}z^{-m-1}
%w^{-n}
=\begin{cases}\delta_{\alpha,\beta}
\,\iota_{z,w}\left(\frac{1}{z-w}\right)
&\quad\text{if}\quad \alpha,\beta\in\Delta^+_r, \\
   0 &\quad\text{otherwise} \end{cases}\label{contract1}\\
\lfloor a_{\alpha}^*a_{\beta}\rfloor
&
=\begin{cases} -\delta_{\alpha,\beta}\,\iota_{z,w}\left(\frac{1}{z-w}
\right)
&\quad\text{if}\quad \alpha,\beta\in\Delta^+_r, \\
-\delta_{\alpha,\beta}\,\delta(w-z)
&\quad\text{otherwise}
\end{cases}
\label{contract2}
\end{align}

As an example of a computation using formal distributions we have the
following

\begin{lem}\label{prelim} For $1\leq i\leq j\leq n$,
$\alpha,\beta\in\Delta^+$,
{\begin{align*}
     [\mathcal H_{i}(z),  a_{\alpha}(w)]
&= -
    (\alpha_i|\alpha)
       a_{\alpha}(z)\delta\left(z-w\right) ,\\
    [\mathcal H_{i}(z) , a_{\alpha}^*(w)]
&=(\alpha_i|\alpha)
       a_{\alpha}^*(z)\delta(z-w),  \\
    [\mathcal H_{i}(z) , \partial_wa_{\alpha}^*(w)]
&=  (\alpha_i|\alpha)
       a_{\alpha}^*(z)\partial_w(z-w) , \\
[\mathcal H_{i}(z),\mathcal H_{ j}(w)]&=
-(\alpha_i|\alpha_j)\left((1-\delta_{i>r}\delta_{j>r})(r+1)+
    \frac{r}{2}\delta_{i,r+1}\delta_{j,r+1}\right)
\partial_w\delta(z-w),
\end{align*}
\begin{align*}
&[\mathcal H_{i}(z) ,:a_\alpha(w)a^*_\beta(w)a^*_\gamma(w):] \\
&\hskip 50pt= (\alpha_i|\beta+\gamma-\alpha)
    :a_{\alpha}(w)a_\beta^*(w)a^*_\gamma(w):\delta(z-w) \\
&\hskip 50pt\quad-\delta_{\alpha\in\Delta^+_r}(\alpha_i|\alpha)
    \left(\delta_{\alpha,\beta} a^*_\gamma(w)+\delta_{\alpha,\gamma}
    a^*_\beta(w)\right)
   \partial_w\delta(z-w),\\ \\
&[\mathcal H_{\alpha_i}(z) ,:a_\alpha(w)a_\beta(w)a^*_\gamma(w):] \\
&\hskip 50pt= (\alpha_i|\gamma-\alpha-\beta)
    :a_{\alpha}(w)a_\beta(w)a^*_\gamma(w):\delta(z-w) \\
&\hskip 50pt\quad-\delta_{\gamma\in\Delta^+_r}(\alpha_i|\gamma)
    \left(\delta_{\gamma,\beta}
a_\gamma(w)+\delta_{\alpha,\gamma}
    a_\beta(w)\right)
    \partial_w\delta(z-w).
\end{align*} }\end{lem}

\begin{proof} %We leave the proof of this Lemma to the reader.
Now by \eqnref{contract1} and \eqnref{contract2} and by Wick's Theorem
\begin{align*}
    \sum_{j}
       :a_{ij}(z) a^*_{ij}(z):
        a_{kl}(w)
&\sim \delta_{ik}
       a_{kl}(z)\lfloor  a^*_{ij}
        a_{kl}\rfloor
    \end{align*}
   and if $\alpha=\alpha_k+\cdots +\alpha_l$, then
\begin{align*}
     \mathcal H_i(z)  a_{kl}(w)
&=\Big(2:a_{ii}a_{ii}^*:+\sum_{j=1}^{i-1}\left(:
          a_{ji}a_{ji}^*:-:a_{j,i-1} a_{j,i-1}^*:\right) \\
        &\quad+\sum_{j=i+1}^{n}\left(:a_{ij}a_{ij}^*:
         - :a_{i+1,j}a_{i+1,j}^*:\right) \Big)
        a_{kl}(w)\\
&\sim -\delta_{1\leq l\leq
r}\left(\delta_{ik}-\delta_{i+1,k}-\delta_{i,l+1}+\delta_{i+1,l+1}\right)
       a_{kl}(z)\iota_{z,w}\left(\frac{1}{z-w}\right) \\
&\qquad-\delta_{r<l}
(\delta_{ik}-\delta_{i+1,k}-\delta_{i,l+1}+\delta_{i+1,l+1})
       a_{kl}(z)\delta(z-w) \\
&\sim \delta_{1\leq l\leq r}(\alpha_i|\alpha)
       a_{kl}(z)\iota_{z,w}\left(\frac{1}{w-
z}\right)-\delta_{r<l}(\alpha_i|\alpha)
       a_{kl}(z)\delta(z-w).
\end{align*}
On the other hand
\begin{align*}
     a_{kl}(w)\mathcal H_i(z)
&\sim \delta_{1\leq l\leq r}(\alpha_i|\alpha)
       a_{kl}(w)\iota_{w,z}\left(\frac{1}{z-w}\right).
\end{align*}
Combining the above operator product expansions we get the first
identity.
A similar computation yields the second identity.
  The third identity comes from differentiating the second with respect
to
$w$.

   On the other hand by Wick's Theorem
\begin{align*}
:a_\nu(z)&a^*_\mu(z)::a_\alpha(w)a^*_\beta(w):
   \\
   &=:a_\alpha(w)a^*_\beta(w)
   a_\nu(z)a^*_\mu(z): +\lfloor a_{\alpha}a_{\mu}^*\rfloor
    :a_\nu(z)a^*_\beta(w):\\
&\quad+
       \lfloor a_{\beta}^*a_{\nu}\rfloor
      :a_\alpha(w)a^*_\mu(z):
      +\lfloor a_{\alpha}a_{\mu}^*\rfloor\lfloor
a_{\beta}^*a_{\nu}\rfloor.
\end{align*}
   Thus
   \begin{align*}
      \mathcal H_{\alpha_i}(z)   \mathcal H_{\alpha_j}(w)
&= \sum_{\alpha,\beta\in\Delta^+}(\alpha_i|\alpha)(\alpha_j|\beta)
       :a_{\alpha}(z) a^*_{\alpha}(z):
       :a_{\beta}(w) a^*_{\beta}(w):\\
&= \sum_{\alpha,\beta\in\Delta^+}(\alpha_i|\alpha)(\alpha_j|\beta)
       :a_{\alpha}(z)a_{\beta}(w) a^*_{\alpha}(z)
a^*_{\beta}(w):\\
&\quad +
\sum_{\beta\in\Delta^+}(\alpha_i|\beta)(\alpha_j|\beta)
        :a_{\beta}(w) a^*_{\beta}(z):\lfloor a_{\beta}^*a_{\beta}\rfloor\\
&\quad +
\sum_{\alpha\in\Delta^+}(\alpha_i|\alpha)(\alpha_j|\alpha)
       :a_{\alpha}(z) a^*_{\alpha}(w) :
       \lfloor a_{\alpha}a_{\alpha}^*\rfloor\\
&\quad+\sum_{\alpha\in\Delta^+_r}(\alpha_i|\alpha)(\alpha_j|\alpha)
\lfloor a_{\alpha}a_{\alpha}^*\rfloor\lfloor
a_{\alpha}^*a_{\alpha}\rfloor,
    \end{align*}
   which can be rewritten as
\begin{align*}
      [\mathcal H_{\alpha_i}(z),\mathcal H_{\alpha_j}(w)]
&=\sum_{\alpha\in\Delta^+_r}(\alpha_i|\alpha)(\alpha_j|\alpha)\left(
\iota_{w,z}
{\frac{1}{(w-z)^2}}-\iota_{z,w}
{\frac{1}{(z-w)^2}}\right) \\
&=-(\alpha_i|\alpha_j)\left((1-\delta_{i>r}\delta_{j>r})(r+1)
     +\frac{r}{2}\delta_{i,r+1}\delta_{j,r+1}\right)
\partial_w\delta(z-w).
    \end{align*}
     This follows from the calculation below for root
system of $\mathfrak{sl}(r+1)$: If $j\leq r$, then
\begin{align*}
    \sum_{\alpha\in\Delta^+_r}(\alpha_j|\alpha)\alpha=
(r+1)\alpha_j
\end{align*}
and
\begin{align*}
    \sum_{\alpha\in\Delta^+_r}(\alpha_{r+1}|\alpha)^2=r.
\end{align*}
    Again by \eqnref{contract1},  \eqnref{contract2} and  Wick's
Theorem
\begin{align*}
:a_{\nu}(z)a^*_{\nu}(z):
:a_\alpha(w)&a^*_\beta(w)a^*_\gamma(w): \\
&=:a_{\nu}(z)a^*_{\nu}(z)a_\alpha(w)
       a^*_\beta(w)a^*_\gamma(w): \\
&\quad +\lfloor a^*_\nu a_\alpha\rfloor
    :a_{\nu}(z)a_\beta^*(w)a^*_\gamma(w): \\
&\quad+
       \lfloor a_\nu a_\beta^*\rfloor
      :a_\alpha(w)a^*_{\nu}(z)a^*_\gamma(w):\\
&\quad+
       \lfloor a_\nu a_\gamma^*\rfloor
      :a_\alpha(w)a^*_\beta(w)a^*_{\nu}(z): \\
   &\quad+
    \left(\lfloor a^*_\nu a_\alpha\rfloor\lfloor a_\nu a_\beta^*\rfloor
a^*_\gamma(w)+\lfloor a^*_\nu a_\alpha\rfloor\lfloor a_\nu
a_\gamma^*\rfloor
    a^*_\beta(w)\right).
\end{align*}

Hence the last identity follows from
\begin{align*}
\mathcal H_{\alpha_i}(z) &:a_\alpha(w)a^*_\beta(w)a^*_\gamma(w):=
       \sum_{\nu\in\Delta^+}(\alpha_i|\nu)
       :a_{\nu}(z) a^*_{\nu}(z):
       :a_\alpha(w)a^*_\beta(w)a^*_\gamma(w):\\
&\sim \sum_{\nu\in\Delta^+}(\alpha_i|\nu)
\Bigg(\Big(\lfloor a^*_\nu a_\alpha\rfloor
    :a_{\nu}(z)a_\beta^*(w)a^*_\gamma(w): \\
&\quad+
       \lfloor a_\nu a_\beta^*\rfloor
      :a_\alpha(w)a^*_{\nu}(z)a^*_\gamma(w): +
       \lfloor a_\nu a_\gamma^*\rfloor
      :a_\alpha(w)a^*_\beta(z)a^*_{\nu}(z)\Big): \\
   &\quad+\left(\lfloor a^*_\nu a_\alpha\rfloor\lfloor a_\nu
a_\beta^*\rfloor
a^*_\gamma(w)+\lfloor a^*_\nu a_\alpha\rfloor\lfloor a_\nu
a_\gamma^*\rfloor
    a^*_\beta(w)\right)  \\
&\sim\Big((\alpha_i|\alpha)\lfloor a^*_\alpha a_\alpha\rfloor
    :a_{\alpha}(z)a_\beta^*(w)a^*_\gamma(w): +
       (\alpha_i|\beta)\lfloor a_\beta a_\beta^*\rfloor
      :a_\alpha(w)a^*_{\beta}(z)a^*_\gamma(w): \\
&\quad+
       (\alpha_i|\gamma)\lfloor a_\gamma a_\gamma^*\rfloor
      :a_\alpha(z)a^*_\beta(z)a^*_\mu(w)\Big): \\
   &\quad+(\alpha_i|\alpha)
    \left(\delta_{\alpha,\beta} a^*_\gamma(w)+\delta_{\alpha,\gamma}
    a^*_\beta(w)\right)
    \lfloor a^*_\alpha a_\alpha\rfloor\lfloor
a_\alpha a_\alpha^*\rfloor.
\end{align*}
\end{proof}

\begin{lem}\label{collection}
\begin{align*}
[a_{ij}(z),a^*_{kl}(w)]&=\delta_{ik}\delta_{jl}\delta(z-w)
\\
[a_{ij}(z)a_{ij}^*(z),a_{ij}(w)a^*_{ij}(w)]&=-\delta_{1\leq i,j\leq r}
\partial_w\delta(z-w) \\
[a_{ij}(z), \partial_w a_{kl}^*(w)]
      &=\delta_{ik}\delta_{jl}\partial_w\delta(z-w)\\
\partial_w a^*_{ij}(w)\delta(z-w)
      &=a^*_{ij}(z)\partial_w\delta(z-w)-a^*_{ij}(w)\partial_w\delta(z-w)
\end{align*}
\end{lem}

The following result collects some other computations
involving the formal distributions that will make future calculations 
less
tedious.
\begin{lem}\label{prelim2}

\begin{align*}
&\sum_{k=i+1,l=j+1}^{n}\left[a_{ik}(z)a_{i+1,k}^*(z),a_{jl}(w)a_{j+1,l}^*
(w)
\right] \tag a\\
&\quad\quad =\Big(\delta_{i,j+1}\sum_{k=j+2}^{n}a_{jk}(z)a_{j+2,k}^*(z)-
    \delta_{j,i+1}\sum_{k=i+2}^{n}a_{ik}(z)a_{i+2,k}^*(z)\Big)\delta(z-
w)\notag
        \\
&\sum_{k=1}^{i-1}\sum_{l=j+1}^{n}\left[a_{k,i-1}(z)a_{ki}^*(z)
        ,a_{jl}(w)a_{j+1,l}^*(w)\right]=
       -\delta_{j,i-1}a_{i-1,i-1}(z)a^*_{ii}(z)\delta(z-w)\tag b \\ \\
&\sum_{l=j+1}^{n}\Big[:a_{ii}^*\left(
        \sum_{k=1}^{i-1} a_{k,i-1}a_{k,i-1}^*-\sum_{k=1}^{i}
        a_{ki}a^*_{ki}\right): ,a_{jl}(w)a_{j+1,l}^*(w)\Big]=0, \tag 
c\\ \\
\end{align*}

\begin{align*}
&\Big[:a_{ii}^*(z)\left(
        \sum_{k=1}^{i-1} a_{k,i-1}(z)a_{k,i-1}^*(z)-\sum_{k=1}^{i}
        a_{ki}(z)a^*_{ki}(z)\right):  ,a_{jj}(w)\Big] \tag d\\
&=-\delta_{ij}\left(:\sum_{k=1}^{i-1}
a_{k,i-1}(z)a_{k,i-1}^*(z)-\sum_{k=1}^{i}
       a_{ki}(z)a^*_{ki}(z):\right)\delta(z-w) \\
&\quad -\left(\delta_{j,i-1}a_{i-1,i-1}(z)a^*_{i,i}(z)-
   \delta_{i,j}:a_{ii}(z)a^*_{ii}(z):\right)\delta(z-w) \\  \\
&\Big[:a_{ii}^*\Big(
        \sum_{k=1}^{i-1} a_{k,i-1}a_{k,i-1}^*\Big):,
      :a_{jj}^*\left(
        \sum_{l=1}^{j-1} a_{l,j-1}a_{l,j-1}^*\right):\Big]\tag e \\
& \quad  =\quad \delta_{j,i-1}:a_{ii}^*(z)a_{i-1,i-1}^*(z)
   \left(
        \sum_{k=1}^{i-2}
     a_{k,i-2}(w)a_{k,i-2}^*(w)\right):\delta(z-w) \\
&    \quad -\delta_{i,j-1}:a_{jj}^*(w)a_{j-1,j-1}^*(w)\left(
        \sum_{l=1}^{j-2} a_{l,j-2}(z)a_{l,j-2}^*(z)\right):
      \delta(z-w) \\
&  \quad -(i-1)\delta_{1\leq i-1\leq
r}\delta_{ij}:a_{ii}^*(z)a_{ii}^*(w):
     \partial_w\delta(z-w) \\ \\
&\Big[:a_{ii}^*\Big(
        \sum_{k=1}^{i-1} a_{k,i-1}a_{k,i-1}^*\Big):,
      :a_{jj}^*\left(
        \sum_{l=1}^{j} a_{l,j}a_{l,j}^*\right):\Big]\tag f \\
&  =\delta_{j,i-1}:a_{ii}^*(z)a_{i-1,i-1}^*(z)\left(
        \sum_{k=1}^{i-1}
      a_{k,i-1}(w)a_{k,i-1}^*(w)\right):\delta(z-w) \\
&  \quad -\delta_{ij}:a_{ii}^*(w)a_{ii}^*(w)\left(
        \sum_{l=1}^{i-1} a_{l,i-1}(z)a_{l,i-1}^*(z)\right):
     \delta(z-w) \\
&   \quad -i\delta_{1\leq i-1\leq r}
     \delta_{j,i-1}:a_{ii}^*(z)a_{i-1,i-1}^*(w):
     \partial_w\delta(z-w) \\ \\
& \Big[:a_{ii}^*\Big(\sum_{k=1}^{i} a_{ki}a_{ki}^*\Big):,
   :a_{jj}^*\left(
        \sum_{l=1}^{j} a_{lj}a_{lj}^*\right):\Big]\tag g  \\
&=-(3+i)\delta_{1\leq
i\leq r}\delta_{ij}:a_{ii}^*(z)a_{ii}^*(w):
     \partial_w\delta(z-w)  \\ \\
&\sum_{l=j+1}^{n}\Big[:a_{ii}^*\left(
        \sum_{k=1}^{i-1} a_{k,i-1}a_{k,i-1}^*-\sum_{k=1}^{i}
        a_{ki}a^*_{ki}\right): ,a_{j+1,l}(w)a_{jl}^*(w)\Big]\tag h\\ \\
&=-\delta_{j,i-1}
     :a^*_{i-1,i}(w)\left(\sum_{k=1}^{i-1}a_{k,i-1}(z)a^*_{k,i-1}(z)
     -\sum_{k=1}^{i}a_{ki}(z)a^*_{ki}(z):\right)\delta(z-w) \\
   &\quad +\delta_{i\leq r}
     \delta_{j,i-1}a_{i-1,i}^*(z)\partial_w\delta(z-w),
\end{align*}

\begin{align*}
&\Big[:a_{ii}^*\left(
        \sum_{k=1}^{i-1} a_{k,i-1}a_{k,i-1}^*-\sum_{k=1}^{i}
        a_{ki}a^*_{ki}\right): ,\sum_{l=1}^{j-1}a_{l,j-1}a_{lj}^*\Big]
\tag i\\
&=\Bigg(:a_{ii}^*\left(\delta_
{i-1,j}\sum_{l=1}^{i-2}a_{l,i-2}a_{l,i-1}^*-
     2\delta_{ij}\sum_{l=1}^{i-1}a_{l,i-1}a_{li}^*
     +\delta_{i,j-1}\sum_{l=1}^{i}a_{l,i}a_{l,i+1}^*\right): \\
& - \delta_{i,j-1}:a^*_{i,i+1}\left(
    \sum_{l=1}^{i-1}a_{l,i-1}a_{l,i-1}^*-
    \sum_{l=1}^{i}a_{li}a_{li}^*\right) :\Bigg)\delta(z-w) \\
&\Big[\sum_{k=1}^
{i-1}a_{k,i-1}a_{ki}^*,\sum_{l=1}^{j-1}a_{l,j-1}a_{lj}^*\Big]
    \tag j\\
&=\Bigg(\delta_{j,i-1}\sum_{l=1}^{i-2}:a_{l,i-2}a_{li}^*:
    -\delta_{i,j-1}\sum_{l=1}^{j-2}:a_{l,j-2}a_{lj}^*:
     \Bigg)\delta(z-w) \\ \\
&\Big[\sum_{k=i+1}^{n}a_{i+1,k}(z)a_{ik}^*(z)
    ,\sum_{l=j+1}^{n} a_{j+1,l}(w)a_{jl}^*(w)\Big]\tag k\\
&=\Bigg(\delta_{i,j-1}\sum_{l=i+2}^{n}:a_{i+2,l}a_{il}^*:-
     \delta_{j,i-1}\sum_{l=j+2}^{n}:a_{j+2,k}a_{jk}^*:
     \Bigg)\delta(z-w) \\ \\
&\Big[\sum_{k=1}^{i-1}a_{k,i-1}a_{ki}^*
    ,\sum_{l=j+1}^{n} a_{j+1,l}
a_{jl}^*\Big]=0.
\tag l\\
\end{align*}
\end{lem}

\section{proof of Theorem \ref{realization}}
We can now check the defining relations.

\begin{lem}[R1]
$$
[\rho(H_i)(z),\rho(H_j)(w)]=(\alpha_i|\alpha_j)\rho(c)\partial_w
\delta(z-w).
$$
\end{lem}
\begin{proof}
We use \lemref{prelim}
in the following calculation:
\begin{align*}[\rho(H_i)&(z),\rho(H_j)(w)]=[\mathcal H_i(z)+b_i(z)
          ,\enspace \mathcal H_j(z)+b_j(z)] \\
&=\Bigg(
-(\alpha_i|\alpha_j)\left((1-\delta_{i>r}\delta_{j>r})(r+1)
    +\frac{r}{2}\delta_{i,r+1}\delta_{j,r+1}\right) \\
&\quad \quad +(\alpha_i|\alpha_j)
    \left( \gamma^2-\delta_{i>r}\delta_{j>r}(r+1)
    +\frac{r}{2}\delta_{i,r+1}\delta_{j,r+1})\right)\Bigg)
     \partial_w\delta(z-w)\\
&=(\alpha_i|\alpha_j)\rho(c)\partial_w \delta(z-w).
\end{align*}

\end{proof}

\begin{lem}[R2]
$$
[\rho(H_i)(z),\rho(E_j)(w)]=(\alpha_i|,\alpha_j)\rho(E_j)(z)\delta(z/w).
$$
\end{lem}

\begin{proof}  We will use \lemref{prelim} repeatedly and the convention
\eqnref{alphakl}:

\begin{align*}
& [\mathcal H_i(z),\rho(E_j)(w)] \\%=[\mathcal H_i(z),\enspace :a_{jj}^*(
%       \sum_{k=1}^{j-1} a_{k,j-1}a_{k,j-1}^*
% -\sum_{k=1}^{j}
%       a_{kj}a^*_{kj}):\\
%&\qquad +\sum_{k=j+1}^n a_{j+1,k}a_{jk}^*
%         -\sum_{k=1}^{j-1}a_{k,j-1}a_{kj}^*
%-a_{jj}^*b_{j}-
%     \left(\delta_{j> r}(r+1)+\delta_{j\leq
%r}(j+1)-\gamma^2\right)\partial a_{jj}^*] \\
&=\sum_{k=1}^{j-1}[\mathcal H_i(z),\enspace :a_{jj}^*
         a_{k,j-1}a_{k,j-1}^*:]-\sum_{k=1}^{j}[\mathcal H_i(z),\enspace
         :a_{jj}^*a_{kj}a^*_{kj}: ]\\
         &\qquad+\sum_{k=j+1}^n [\mathcal H_i(z),\enspace
a_{j+1,k}a_{jk}^*]
          -\sum_{k=1}^{j-1}[\mathcal H_i(z),\enspace a_{k,j-1}a_{kj}^*] \\
&\qquad-[\mathcal H_i(z),\enspace a_{jj}^*]b_{j}-
      \left(\delta_{j> r}(r+1)+\delta_{j\leq
r}(j+1)-\gamma^2\right)[\mathcal H_i(z),\enspace \partial_w
a_{jj}^*(w)] \\ \\
&=\sum_{k=1}^{j-1}\left((\alpha_i|\alpha_j):a_{jj}^*
         a_{k,j-1}a_{k,j-1}^*:\delta(z-w)-\delta_{j-1\leq  r}
  (\alpha_i|\alpha_{k,j-1})a_{jj}^*(w)\partial_w\delta(z-w)\right)\\
&\quad-\sum_{k=1}^{j}\left((\alpha_i|\alpha_j)
         :a_{jj}^*a_{kj}a^*_{kj}:\delta(z-w) -\delta_{ j\leq r}
   (\alpha_i|\alpha_{kj})
   (\delta_{jk}a_{kk}^*(w)+a_{jj}^*(w))\partial_w\delta(z-w)\right)\\
&\quad+(\alpha_i|\alpha_{j})\sum_{k=j+1}^n
  a_{j+1,k}(z)a_{jk}^*(w)\delta(z-w)
      -(\alpha_i|\alpha_j)\sum_{k=1}^{j-1}
  a_{k,j-1}(z)a_{kj}^*(w)\delta(z-w) \\
&\quad-(\alpha_i|\alpha_j)a_{jj}^*(z)b_{j}(w)\delta(z-w) \\
&\quad-(\alpha_i|\alpha_j)
      \left(\delta_{j> r}(r+1)+\delta_{j\leq
r}(j+1)-\gamma^2\right)a_{jj}^*(z)\partial_w\delta(z-w)
\end{align*}
\begin{align*}
&=(\alpha_i|\alpha_j)\Bigg(:a_{jj}^*(\sum_{k=1}^{j-1}
         a_{k,j-1}a_{k,j-1}^*-\sum_{k=1}^{j}a_{kj}a^*_{kj}):\\
&\quad +\sum_{k=j+1}^n
a_{j+1,k}a_{jk}^*
     -\sum_{k=1}^{j-1} a_{k,j-1}a_{kj}^* 
-a_{jj}^*b_{j}\Bigg)\delta(z-w) \\
&\qquad+\sum_{k=1}^{j}\delta_{1\leq j\leq
r}(\alpha_i|\alpha_k+\cdots
         +\alpha_{j})(\delta_{jk}a_{kk}^*(w)+a_{jj}^*(w))\partial_w\delta(z-
w)\\
&\quad-\sum_{k=1}^{j-1}\left(\delta_{1\leq j-1\leq
r}(\alpha_i|\alpha_k+\cdots
    +\alpha_{j-1})a_{jj}^*(w)\partial_w\delta(z-w)\right)\\
&\qquad-(\alpha_i|\alpha_j)
      (\delta_{j> r}(r+1)+\delta_{j\leq
r}(j+1)-\gamma^2)a_{jj}^*(z)\partial_w\delta(z-w).
\end{align*}

The last term of $\rho(H_i)(z)$ gives us
\begin{align*}
[b_i(z),\enspace
\rho(E_j)(w)]&=-\mathfrak{B}_{ij}a_{jj}(w)\partial_w\delta(z-w).
\end{align*}

There are three cases to consider:

\noindent Case I.  $j\leq r$:  Then

\begin{align*}
& \sum_{k=1}^{j}(\alpha_i|\alpha_k+\cdots
 +\alpha_{j})(\delta_{jk}a_{kk}^*(w)+a_{jj}^*(w))\partial_w\delta(z-w)\\
& -\sum_{k=1}^{j-1}(\alpha_i|\alpha_k+\cdots
    +\alpha_{j-1})a_{jj}^*(w)\partial_w\delta(z-w)\\
& -(\alpha_i|\alpha_j)
      \left(\delta_{j> r}(r+1)+\delta_{j\leq
r}(j+1)-\gamma^2\right)a_{jj}^*(z)\partial_w\delta(z-w)\\
& -(\alpha_i|\alpha_j)(\gamma^2
     -\delta_{i>r}\delta_{j>r}(r+1)
     +\frac{r}{2} \delta_{i,r+1}\delta_{j,r+1})a_{jj}(w)\partial_w\delta
(z-w)  \\
& =(j+1)(\alpha_i|\alpha_j)a_{jj}^*(w)\partial_w\delta(z-w)\\
& -(\alpha_i|\alpha_j)
      (j+1-\gamma^2)a_{jj}^*(z)\partial_w\delta(z-w)\\
& -(\alpha_i|\alpha_j)\gamma^2a_{jj}(w)\partial_w\delta(z-w)  \\
& =-(\alpha_i|\alpha_j)(j+1-\gamma^2)\partial_wa_{jj}(w)\delta(z-w)
\end{align*}
by \lemref{collection}.

\noindent Case II: $j=r+1$:
\begin{align*}
&-\sum_{k=1}^{r}\left((\alpha_i|\alpha_k+\cdots
    +\alpha_{r})a_{r+1,r+1}^*(w)\partial_w\delta(z-w)\right)\\
&-(\alpha_i|\alpha_{r+1})
      \left((r+1)-\gamma^2\right)a_{r+1,r+1}^*(z)\partial_w\delta(z-w)\\
&-(\alpha_i|\alpha_{r+1})(\gamma^2
     -\delta_{i>r}(r+1)
     +\frac{r}{2}\delta_{i,r+1})a_{r+1,r+1}(w)\partial_w\delta(z-w) \\
&=-(\alpha_i|\alpha_{r+1})
      \left((r+1)-\gamma^2\right)(\partial_wa_{r+1,r+1}^*(w))\delta(z-w)\\
\end{align*}
which follows from
\begin{align*}
&-\sum_{k=1}^{r}(\alpha_i|\alpha_k+\cdots
    +\alpha_{r})+(\alpha_i|\alpha_{r+1})(
     \delta_{i>r}(r+1)
     -\frac{r}{2}\delta_{i,r+1})\\
&=\begin{cases}-2 &\quad\text{if}\quad 1=i=r \\
         0&\quad\text{if}\quad 1\leq i< r \\
         -(r+1)&\quad\text{if}\quad 1< i= r \\
         2(r+1)&\quad\text{if}\quad 1\leq i= r+1 \\
         (\alpha_i|\alpha_{r+1})(r+1)&\quad\text{if}\quad i> r+1
     \end{cases} \\
&=(\alpha_i|\alpha_{r+1})(r+1)\\
\end{align*}

Cases III:  $j>r+1$:
\begin{align*}
&-(\alpha_i|\alpha_j)
      \left(r+1-\gamma^2\right)a_{jj}^*(z)\partial_w\delta(z-w)\\
&-\mathfrak{B}_{ij}a_{jj}(w)\partial_w\delta(z-w) \\
&=-(\alpha_i|\alpha_j)
      \left(r+1-\gamma^2\right)a_{jj}^*(z)\partial_w\delta(z-w)  \\
&\quad +(\alpha_i|\alpha_j)(\gamma^2
     -\delta_{i>r}(r+1))a_{jj}(w)\partial_w\delta(z-w) \\
&=-(\alpha_i|\alpha_j)
      \left(r+1-\gamma^2\right)(\partial_wa_{jj}^*(w))\delta(z-w)
\end{align*}
by \lemref{collection} and the fact that $(\alpha_i|\alpha_j)=0$ for
$i\leq r<r+1<j$.

Putting these computations together we get
\begin{align*}
[\rho( H_i)(z)&,\rho(E_j)(w)]=(\alpha_i|\alpha_j)\rho(E_j)(w)\delta(z-w).
\end{align*}

\end{proof}

\begin{lem}[R3]
$$
[\rho(H_i)(z),\rho(F_j)(w)]=-(\alpha_i|\alpha_j)\rho(F_j)(z)\delta(z-w).
$$
\end{lem}

\begin{proof}  The proof follows from \lemref{prelim} :
\begin{align*}
[\rho(H_i)(z)&,\rho(F_j)(w)]=[\mathcal
  H_i(z),a_{jj}(w)+\sum_{k=j+1}^{n}a_{jk}(w)a_{j+1,k}^*(w)] \\
                 &=-(\alpha_i|\alpha_j)a_{jj}(w)\delta(z-
w) \\
&\quad +\sum_{k=j+1}^{n}[\mathcal
H_i(z),a_{jk}(w)]a_{j+1,k}^*(w)+\sum_{k=j+1}^{n}a_{jk}(w)[\mathcal
H_i(z),a_{j+1,k}^*(w)]  \\
                 &=\Big
(-(\alpha_i|\alpha_j)a_{jj}(w)-\sum_{k=j+1}^{n}(\alpha_i|\alpha_j+\cdots
+\alpha_k)a_{jk}(z)a_{j+1,k}^*(w)   \\
    &\hskip 50pt+\sum_{k=j+1}^{n}(\alpha_i|\alpha_{j+1}+\cdots
+\alpha_k)a_{jk}(w)a_{j+1,k}^*(z)\Big)\delta(z-w) \\
                 &=-(\alpha_i|\alpha_j)\rho(F_j)(z)\delta(z-w)
\end{align*}
\end{proof}

\begin{lem}[R4]
$$
[\rho(E_i)(z),\rho(F_j)(w)]
=\delta_{i,j}(\rho(H_i)(z))\delta(z-w)+\rho(c)\partial_w\delta(z-w))
$$
\end{lem}

\begin{proof}
First we take $i=j$.  Now for the convenience of the reader
we recall that $\rho(E_i)(z)$ is equal to
\begin{align*} &:a_{ii}^*\left(
        \sum_{k=1}^{i-1} a_{k,i-1}a_{k,i-1}^*-\sum_{k=1}^{i}
        a_{ki}a^*_{ki}\right): +\sum_{k=i+1}^n a_{i+1,k}a_{ik}^*
          -\sum_{k=1}^{i-1}a_{k,i-1}a_{ki}^* \\
      &\quad-a_{ii}^*b_{i}
      -\left(\delta_{i>r}(r+1) +\delta_{i\leq
r}(i+1)-\gamma^2\right)\partial
a_{ii}^*\\
\end{align*}
and thus the first summand of
$\rho(F_i)(w)=a_{ii}+\sum_{l=i+1}^{n}a_{il}a_{i+1,l}^*$ brackets with
$\rho(E_i)(z)$ to give us (by \lemref{prelim2} (d) and
\lemref{collection})

\begin{align*}
&\quad \Big(2:a_{ii}(z)a^*_{ii}(z):-:
        \sum_{k=1}^{i-1}\left( a_{k,i-1}a_{k,i-1}^*-
        a_{ki}a^*_{ki}\right):+ b_{i}(z)\Big)
    \delta(z-w)  \\
      &\quad+\left(\delta_{i>r}(r+1) +\delta_{i\leq
r}(i+1)-\gamma^2\right)\partial_z
\delta(z-w).
\\
\end{align*}
The second summation in $\rho(F_i)(w)$ contributes

\begin{align*}
\sum_{l=i+1}^{n}\Big[\rho(E_i)(z),
&\enspace a_{il}(w)a_{i+1,l}^*(w)]=
\sum_{l=i+1}^{n}\Big[\sum_{k=i+1}^n a_{i+1,k}(z)a_{ik}^*(z),
         a_{il}(w)a_{i+1,l}^*(w)\Big] \\
       &=\sum_{l=i+1}^{n}\Big(a_{il}(z)a_{il}^*(z)-
        a_{i+1,l}(z)a_{i+1,l}^*(z)\Big)
        \delta(z-w) \\
       &\quad -\delta_{i+1\leq r}(r-i  )\partial_w\delta(z-w).
\end{align*}
Adding these two summations up, we arrive at the desired result.

Now consider the case $|i-j|\geq 1$. Then $\rho(F_j)(w)$ is
$
a_{jj}+\sum_{l=j+1}^{n}a_{jl}a_{j+1,l}^*.
$
First we have

\begin{align*} [E_i(z),a_{jj}(w)]%=\Big[:a_{ii}^*\left(
%       \sum_{k=1}^{i-1} a_{k,i-1}a_{k,i-1}^*-\sum_{k=1}^{i}
%       a_{ki}a^*_{ki}\right): +\sum_{k=i+1}^n a_{i+1,k}a_{ik}^*
%         -\sum_{k=1}^{i-1}a_{k,i-1}a_{ki}^* \\
&=\Big[:a_{ii}^*\left(
        \sum_{k=1}^{i-1} a_{k,i-1}a_{k,i-1}^*-\sum_{k=1}^{i}
        a_{ki}a^*_{ki}\right): , a_{jj}(w)\Big]\\
    &=-\delta_{j,i-1}a_{i-1,i-1}(z)a^*_{i,i}(z)\delta(z-w)
\end{align*}
by \lemref{prelim2} (d).
Next we have
\begin{align*}
&[E_i(z),\sum_{l=j+1}^{n}a_{jl}(w)a_{j+1,l}^*(w)]\\
&=\Big[:a_{ii}^*\left(
        \sum_{k=1}^{i-1} a_{k,i-1}a_{k,i-1}^*-\sum_{k=1}^{i}
        a_{ki}a^*_{ki}\right): +\sum_{k=i+1}^n a_{i+1,k}a_{ik}^*
          -\sum_{k=1}^{i-1}a_{k,i-1}a_{ki}^* \\
      &\quad-a_{ii}^*b_{i}
      +\left(\gamma^2 -\delta_{i+1\leq
r}(i+1)\right)\partial
a_{ii}^*,\sum_{l=j+1}^{n}a_{jl}(w)a_{j+1,l}^*(w)\Big]\\
&=\Big[:a_{ii}^*\left(
        \sum_{k=1}^{i-1} a_{k,i-1}a_{k,i-1}^*-\sum_{k=1}^{i}
        a_{ki}a^*_{ki}\right):
          -\sum_{k=1}^{i-1}a_{k,i-1}a_{ki}^*
,\sum_{l=j+1}^{n}a_{jl}(w)a_{j+1,l}^*(w)\Big]\\
&\hskip 100pt  \text{as $i\neq j$,}\\
    &=\delta_{j,i-1}a_{i-1,i-1}(z)a^*_{ii}(z)\delta(z-w) \\
&\hskip 100pt  \text{by
\lemref{prelim2} (b) and (c).}\\
\end{align*}

Adding up the last two calculations finishes the proof of this lemma.

\end{proof}

We are now left with the Serre relations:

\begin{lem}[R5/R6]

\begin{align*}
[\rho(F_i)(z),\rho(F_j)(w)]&=[\rho(E_i)(z),\rho(E_j)(w)]=0
    \quad\text{if}\quad
      (\alpha_i|\alpha_j)\neq -1 \\
   [\rho(F_i)(z_1),\rho(F_i)(z_2),\rho(F_j)(w)]
     &=[\rho(E_i)(z_1),\rho(E_i)(z_2),\rho(E_j)(w)]=0, \\
    &\quad\text{if}\quad  (\alpha_i|\alpha_j)= -1.
\end{align*}
\end{lem}

\begin{proof}
Let us check the relations for $\rho(F_i)$. (The Serre relations were
already
known to hold true for the $F_i$, see B. Feigin and E. Frenkel$^7$,
but we  provide a
proof as some of the calculations will be used in future work.)
%A straight forward  calculation shows
%
%\begin{align*}\Big[a_{ii}(z),\enspace
%   &a_{jj}+\sum_{l=j+1}^{n}a_{jl}a_{j+1,l}^*
%   \Big] =\delta_{i,j+1}a_{j,j+1}(w)\delta(z-w). \\
%\end{align*}
%  Moreover
%
%\begin{align*}\Big[\sum_{k=i+1}^{n}a_{ik}&(z)a_{i+1,k}^*(z),\enspace
%   a_{jj}(w)+\sum_{l=j+1}^{n}a_{jl}(w)a_{j+1,l}^*(w)
%   \Big] \\
%&=\Big(\delta_{i,j+1}\sum_{k=j+2}^{n}a_{jk}a_{j+2,k}^*-
%   \delta_{j,i+1}\sum_{k=i+2}^{n}a_{ik}a_{i+2,k}^*-
%\delta_{j,i+1}a_{i,i+1}\Big)\delta(z-w)
%\end{align*}
By \lemref{prelim2} (a)

\begin{align}[\rho(F_i)(z),&\rho(F_j)(w)]
=\left(\delta_{i,j+1}a_{i-1,i}-\delta_{j,i+1}a_{i,i+1}
    \right)\delta(z-w)\\
        &\quad
+\left(\delta_{i,j+1}\sum_{k=i+1}^{n}a_{i-1,k}a_{i+1,k}^*-
    \delta_{j,i+1}\sum_{k=i+2}^{n}a_{ik}a_{i+2,k}^*
\right)\delta(z-w).\notag
\end{align}
Note the above is zero if $|i-j|\neq 1$ which is precisely when
$(\alpha_i|\alpha_j)\neq -1$.   As the first index in
$a_{kl}$ (resp.
$a^*_{kl}$) in
$\rho(F_i)(z)$ is
$i$ (resp. $i+1$) we also get
$$
[\rho(F_i)(z_1),\rho(F_i)(z_1),\rho(F_j)(w)]=0.
$$
This completes the proof of the relations R5 and R6 for
$\rho(F_i)(z)$.

Now we break up $\rho(E_i)(z)$ into three summands

\begin{align*}  &\rho(E^1_i)(z):=   :a_{ii}^*\left(
        \sum_{k=1}^{i-1} a_{k,i-1}a_{k,i-1}^*-\sum_{k=1}^{i}
        a_{ki}a^*_{ki}\right):  \\
      &\rho(E^2_i)(z):=\sum_{k=i+1}^n a_{i+1,k}a_{ik}^*
          -\sum_{k=1}^{i-1}a_{k,i-1}a_{ki}^* \\
      &\rho(E^3_i)(z):= -a_{ii}^*b_{i}
      -\left(\delta_{i>r}(r+1)+\delta_{i\leq
r}(i+1)-\gamma^2\right)\partial a_{ii}^*  .
\end{align*}
By \lemref{prelim2} (e), (f) and (g) we have
\begin{align*}
&\Big[\rho(E^1_i)(z),  \rho(E^1_j)(w)\Big]\\
&=\delta_{j,i-1}:a_{ii}^*(z)a_{i-1,i-1}^*(z)\left(
        \sum_{k=1}^{i-2}
     a_{k,i-2}(w)a_{k,i-2}^*(w)\right):\delta(z-w) \\
&-\delta_{i,j-1}:a_{jj}^*(w)a_{j-1,j-1}^*(w)\left(
        \sum_{l=1}^{j-2} a_{l,j-2}(z)a_{l,j-2}^*(z)\right)
      \delta(z-w) \\
&-\delta_{j,i-1}:a_{ii}^*(z)a_{i-1,i-1}^*(z)\left(
        \sum_{k=1}^{i-1}
      a_{k,i-1}(w)a_{k,i-1}^*(w)\right):\delta(z-w) \\
&+\delta_{i,j-1}:a_{jj}^*(z)a_{j-1,j-1}^*(z)\left(
        \sum_{k=1}^{j-1}
      a_{k,j-1}(w)a_{k,j-1}^*(w)\right):\delta(z-w) \\
&-(i-1)\delta_{1\leq i-1\leq r}
   \delta_{ij}:a_{ii}^*(z)a_{ii}^*(w):
     \partial_w\delta(z-w) \\
&+i\delta_{1\leq i-1\leq r}
    \delta_{j,i-1}:a_{ii}^*(z)a_{i-1,i-1}^*(w):
     \partial_w\delta(z-w) \\
&-j\delta_{1\leq j-1\leq r}
     \delta_{i,j-1}:a_{jj}^*(w)a_{j-1,j-1}^*(z):
     \partial_z\delta(z-w) \\
&-(3+i)\delta_{1\leq i\leq r}\delta_{ij}:a_{ii}^*(z)a_{ii}^*(w):
     \partial_w\delta(z-w).
\end{align*}

By \lemref{prelim2} (h) and (i),
\begin{align*}
&\Big[\rho(E^1_i)(z),\rho(E^2_j)(w)\Big]+\Big[\rho(E^2_i)
(z),
\rho(E^1_j)(w)\Big] \\
% &=-\delta_{j,i-1}
%    :a^*_{i-1,i}(w)\left(\sum_{k=1}^{i-1}a_{k,i-1}(z)a^*_{k,i-1}(z)
%    -\sum_{k=1}^{i}a_{ki}(z)a^*_{ki}(z):\right)\delta(z-w) \\
%&\quad+\delta_{i,j-1}
%    :a^*_{j-1,j}(z)\left(\sum_{k=1}^{j-1}a_{k,j-1}(w)a^*_{k,j-1}(w)
%    -\sum_{k=1}^{j}a_{kj}a^*_{kj}:\right)\delta(z-w) \\
%&\quad-\Bigg(:a_{ii}^*\left(\delta_{i-1,j}\sum_{l=1}^{i-2}a_{l,i-2}a_{l,i-1}
%^*
%    +\delta_{i,j-1}\sum_{l=1}^{i}a_{l,i}a_{l,i+1}^*\right) \\
%&\quad - \delta_{i,j-1}a^*_{i,i+1}\left(
%   \sum_{l=1}^{i-1}a_{l,i-1}a_{l,i-1}^*-
%   \sum_{l=1}^{i}a_{li}a_{li}^*\right) :\Bigg)\delta(z-w) \\
%&\quad+\Bigg(:a_{jj}^*\left(\delta_
%{j-1,i}\sum_{l=1}^{j-2}a_{l,j-2}a_{l,j-1}^*
%    +\delta_{j,i-1}\sum_{l=1}^{j}a_{l,j}a_{l,j+1}^*\right) \\
%& \quad- \delta_{j,i-1}a^*_{j,j+1}\left(
%   \sum_{l=1}^{j-1}a_{l,j-1}a_{l,j-1}^*-
%   \sum_{l=1}^{j}a_{lj}a_{lj}^*\right) :\Bigg)\delta(z-w) \\
%&\quad +\delta_{i\leq r}
%    \delta_{j,i-1}a_{i-1,i}^*(z)\partial_w\delta(z-w) \\
%&\quad - \delta_{j\leq r}
%   \delta_{i,j-1}a_{j-1,j}^*(w)\partial_z\delta(z-w). \\
%
&=-\delta_{j,i-1}
     :a^*_{i-1,i}\left(\sum_{k=1}^{i-2}a_{k,i-2}a^*_{k,i-2}
     -\sum_{k=1}^{i}a_{ki}a^*_{ki}:\right)\delta(z-w) \\
&\quad+\delta_{i,j-1}
     :a^*_{i,i+1}\left(
   \sum_{k=1}^{i-1}a_{k,i-1}a^*_{k,i-1}
     -\sum_{k=1}^{i+1}a_{k,i+1}a^*_{k,i+1}:\right)\delta(z-w) \\
&\quad-:a_{ii}^*\left(\delta_{i-1,j}\sum_{l=1}^{i-2}a_{l,i-2}a_{l,i-1}^*
     +\delta_{i,j-1}\sum_{l=1}^{i}a_{l,i}a_{l,i+1}^*\right)
     \delta(z-w) \\
&\quad+:a_{jj}^*\left(\delta_{j-1,i}\sum_{l=1}^{j-2}a_{l,j-2}a_{l,j-1}^*
     +\delta_{j,i-1}\sum_{l=1}^{j}a_{l,j}a_{l,j+1}^*\right)
     \delta(z-w) \\
&\quad +\delta_{i\leq r}
     \delta_{j,i-1}a_{i-1,i}^*(z)\partial_w\delta(z-w) \\
&\quad - \delta_{j\leq r}
    \delta_{i,j-1}a_{j-1,j}^*(w)\partial_z\delta(z-w). \\
\end{align*}
Similarly
% \begin{align*}
%&\Big[\rho(E^1_i)(z)\rho(E^3_j)(w) \Big]\\
% &=\Big[:a_{ii}^*\left(
%       \sum_{k=1}^{i-1} a_{k,i-1}a_{k,i-1}^*-\sum_{k=1}^{i}
%       a_{ki}a^*_{ki}\right): ,
%       -a_{jj}^*b_{j}
%     -\left(\delta_{j>r}(r+1)+\delta_{j\leq
%r}(j+1)-\gamma^2\right)\partial
%a_{jj}^* \Big]\\
%&=\left(-\delta_{i-1,j}a_{ii}^*
%        a_{i-1,i-1}^*+\delta_{ij}
%       :a_{ii}^*a^*_{ii}:\right)b_{j}\delta(z-w)\\
%& +\left(\delta_{j>r}(r+1)+
%    \delta_{j\leq r}(j+1)-\gamma^2\right)\left
%        a_{i-1,i-1}^*(z)+\delta_{ij}
%       :a_{ii}^*(z)a^*_{ii}(z):\right)\partial_w\delta(z-w)\\
%
%\end{align*}
%so that
\begin{align*}
&\Big[\rho(E^1_i)(z)\rho(E^3_j)(w)\Big]
   +\Big[\rho(E^3_i)(z),\rho(E^1_j)(w)\Big] \\
  &=\left(-\delta_{i-1,j}a_{ii}^*
         a_{i-1,i-1}^*b_j+\delta_{j-1,i}a_{jj}^*
         a_{j-1,j-1}^*b_i\right)\delta(z-w)\\
& +\left(\delta_{j>r}(r+1)+
      \delta_{j\leq r}(j+1)-\gamma^2\right)\left
(-\delta_{i-1,j}a_{ii}^*(z)
         a_{i-1,i-1}^*(z)+\delta_{ij}
        a_{ii}^*(z)a^*_{ii}(z)\right)\partial_w\delta(z-w)\\
& +\left(\delta_{i>r}(r+1)+\delta_{i\leq r}(i+1)-\gamma^2\right)
        \left(\delta_{j-1,i}a_{jj}^*(w)
         a_{j-1,j-1}^*(w)-\delta_{ij}
        a_{jj}^*(w)a^*_{jj}(w)\right)\partial_z\delta(z-w).
\end{align*}

By \lemref{prelim2} (j), (k) and (l) we have

\begin{align*}\Big[\rho(E^2_i)(z), \rho(E^2_j)(w)\Big] \\
&=\Bigg(\delta_{j,i-1}\sum_{l=1}^{i-2}:a_{l,i-2}a_{li}^*:
    -\delta_{i,j-1}\sum_{l=1}^{j-2}:a_{l,j-2}a_{lj}^*:
     \Bigg)\delta(z-w) \\
&+\Bigg(\delta_{i,j-1}\sum_{l=i+2}^{n}:a_{i+2,l}a_{il}^*:-
     \delta_{j,i-1}\sum_{l=j+2}^{n}:a_{j+2,k}a_{jk}^*:
     \Bigg)\delta(z-w) \\
%&+\Big(\delta_{j,i+1}\sum_{l=1}^{j-2}
%   a_{l+1,j-1}a_{l,j}^*+\delta_{i,j+1}\sum_{l=1}^{i-2}
%   a_{l+1,i-1}a_{l,i}^*\Big)\delta(z-w)
\end{align*}

Next we have

%
%\begin{align*}
%&\Big[\rho(E^2_i)(z),
%\rho(E^3_j)(w)\Big]\\
%   &=\Big[\sum_{k=i+1}^n a_{i+1,k}a_{ik}^*
%         -\sum_{k=1}^{i-1}a_{k,i-1}a_{ki}^*,
%       -a_{jj}^*b_{j}
%     -\left(\delta_{j>r}(r+1)+
%    \delta_{j\leq r}(j+1)-\gamma^2\right)\partial a_{jj}^* \Big] \\
%   &=\left(-\delta_{j,i+1}a_{i,i+1}^*b_{i+1}
%         +\delta_{j,i-1}a_{i-1,i}^*b_{i-1}\right)\delta(z-w) \\
%   &\quad
%   -\left(\delta_{j>r}(r+1)+\delta_{j\leq r}(j+1)-\gamma^2\right)
%   \left(\delta_{j,i+1}a_{i,i+1}^*(z)
%         -
%     \delta_{j,i-1}a_{i-1,i}^*(z)\right)\partial_w\delta(z-w)
%\end{align*}
%so that
\begin{align*}
&\Big[\rho(E^2_i)(z),
\rho(E^3_j)(w)\Big]+\Big[\rho(E^3_i)(z),
\rho(E^2_j)(w)\Big]\\
    &=\left(-\delta_{j,i+1}a_{i,i+1}^*b_{i+1}
          +\delta_{j,i-1}a_{i-1,i}^*b_{i-1}\right)\delta(z-w) \\
    &\quad-
    \left(\delta_{j>r}(r+1)+\delta_{j\leq r}(j+1)-\gamma^2\right)
    \left(\delta_{j,i+1}a_{i,i+1}^*(z)
         -
      \delta_{j,i-1}a_{i-1,i}^*(z)\right)\partial_w\delta(z-w)  \\
    &-\left(-\delta_{i,j+1}a_{j,j+1}^*b_{j+1}
          +\delta_{i,j-1}a_{j-1,j}^*b_{j-1}\right)\delta(z-w) \\
    &\quad+
    \left(\delta_{i>r}(r+1)+\delta_{i\leq r}(i+1)-\gamma^2\right)
    \left(\delta_{i,j+1}a_{j,j+1}^*(w)
          -
      \delta_{i,j-1}a_{j-1,j}^*(w)\right)\partial_z\delta(z-w).
\end{align*}
while
\begin{align*}
& \Big[\rho(E^3_i)(z),\enspace \rho(E^3_j)(w)\Big]\\
&= \Big[ -a_{ii}^*b_{i}
      -\left(\delta_{i>r}(r+1)+\delta_{i\leq
r}(i+1)-\gamma^2\right)\partial
    a_{ii}^* , \\
&\hskip 75pt -a_{jj}^*b_{j}
      -\left(\delta_{j>r}(r+1)+\delta_{j\leq
r}(j+1)-\gamma^2\right)\partial
    a_{jj}^* \Big]\\
&=a_{ii}^*(z)a_{jj}^*(w)
     \mathfrak{B}_{ij}
\partial_w\delta(z-w).
\end{align*}

Now we observe that every reduction above is zero if $|i-j|\not\in
\{0,1\}$ i.e. if $(\alpha_i|\alpha_j)=0$.  Thus

\begin{align*}
&\Big[\rho(E_i)(z), \enspace \rho(E_j)(w)\Big] =0,\quad \text{if}\quad
(\alpha_i|\alpha_j)=0.
\end{align*}

When $i=j$, $\Big[\rho(E_i)(z), \enspace \rho(E_j)(w)\Big]$ reduces, by
\lemref{collection}, to

\begin{align*}
&\Big[\rho(E_i)(z),\enspace \rho(E_i)(w)\Big] \\
&=-(i-1)\delta_{1\leq i-1\leq r}:a_{ii}^*(z)a_{ii}^*(w):
     \partial_w\delta(z-w) \\
&-(3+i)\delta_{1\leq i\leq r}:a_{ii}^*(z)a_{ii}^*(w):
     \partial_w\delta(z-w) \\
&+\left(\delta_{i>r}(r+1)+
      \delta_{i\leq r}(i+1)-\gamma^2\right)
        a_{ii}^*(z)a^*_{ii}(z)\partial_w\delta(z-w)\\
&-\left(\delta_{i>r}(r+1)+\delta_{i\leq r}(i+1)-\gamma^2\right)
        a_{ii}^*(w)a^*_{ii}(w)\partial_z\delta(z-w) \\
&+2a_{ii}^*(z)a_{ii}^*(w)
      \left( \gamma^2+(\delta_{1\leq i\leq r}-1)(r+1)+
      \delta_{i,r+1}\frac{r}{2}\right)
     \partial_w\delta(z-w) \\
&=-\left((i-1)\delta_{1\leq i-1\leq
r}+(3+i)\delta_{1\leq i\leq r}
-\delta_{i,r+1}r\right):a_{ii}^*(z)a_{ii}^*(w):
     \partial_w\delta(z-w) \\
&+\delta_{i\leq r}(i+1)\left(
        a_{ii}^*(z)a^*_{ii}(z)+
        a_{ii}^*(w)a^*_{ii}(w)\right)\partial_z\delta(z-w) \\
&+(r+1)\delta_{i>r}\left(a_{ii}^*(z)a^*_{ii}(z)+
        a_{ii}^*(w)a^*_{ii}(w)-2a_{ii}^*(z)a_{jj}^*(w)\right)
     \partial_w\delta(z-w) \\
&+\gamma^2
      \left(2a_{ii}^*(z)a_{ii}^*(w)
    -a_{ii}^*(z)a^*_{ii}(z)-a_{ii}^*(w)a^*_{ii}(w)\right)
\partial_w\delta(z-w)=0.
\end{align*}
This proves the result for $i\neq j\pm 1$.

If $i=j+1$ then we get

\begin{align*}
&\Big[\rho(E_i)(z),\enspace \rho(E_{i-1})(w)\Big] \\
&=:a_{ii}^*(z)a_{i-1,i-1}^*(z)\left(
        \sum_{k=1}^{i-2}
     a_{k,i-2}(w)a_{k,i-2}^*(w)
     -\sum_{k=1}^{i-1}
      a_{k,i-1}(w)a_{k,i-1}^*(w)\right):\delta(z-w) \\
&-
     :a^*_{i-1,i}\left(\sum_{k=1}^{i-2}a_{k,i-2}a^*_{k,i-2}
     -\sum_{k=1}^{i}a_{ki}a^*_{ki}:\right)\delta(z-w) \\
&-:a_{ii}^*\left(\sum_{l=1}^{i-2}a_{l,i-2}a_{l,i-1}^*
     \right)
     \delta(z-w)
+:a_{i-1,i-1}^*\left(\sum_{l=1}^{i-1}a_{l,i-1}a_{li}^*\right)
     \delta(z-w) \\
&+\Bigg(\sum_{l=1}^{i-2}:a_{l,i-2}a_{li}^*:
     \Bigg)\delta(z-w) -\Bigg(\sum_{l=i+1}^{n}:a_{i+1,k}a_{i-1,k}^*:
     \Bigg)\delta(z-w) \\ \\
&-a_{ii}^*a_{i-1,i-1}^*b_{i-1}\delta(z-w)
     +a_{i-1,i}^*(b_{i-1}+b_i)\delta(z-w) \\
%&+\delta_{i\leq r}a_{i-1,i}^*(z)\partial_w\delta(z-w) \\
%&+\left(\delta_{i-1>r}(r+1)+\delta_{i-1\leq r}i-\gamma^2\right)
%     a_{i-1,i}^*(z)\partial_w\delta(z-w)  \\
%&-\left(\delta_{i>r}(r+1)+\delta_{i\leq r}(i+1)-\gamma^2\right)
%     a_{i-1,i}^*(w)\partial_w\delta(z-w)\\
%
%&-\left(\delta_{i-1>r}(r+1)+\delta_{i-1\leq r}i
%   -\gamma^2\right)a_{ii}^*(z)
%       a_{i-1,i-1}^*(z)\partial_w\delta(z-w)\\
%&+i\delta_{1\leq i-1\leq r}
%   :a_{ii}^*(z)a_{i-1,i-1}^*(w):
%   \partial_w\delta(z-w) \\
%&-a_{ii}^*(z)a_{i-1,i-1}^*(w)
%    \left( \gamma^2-\delta_{i-1>r}(r+1)\right)
%   \partial_w\delta(z-w). \\ \\
&+\left(\delta_{i>r}(r+1)+\delta_{i\leq r}(i+1)-\gamma^2\right)
       \partial_wa_{i-1,i}^*(w)\delta(z-w)\\
&+a_{ii}^*(z)\partial_wa_{i-1,i-1}^*(w)
      \left( \gamma^2-\delta_{i-1\leq r}i-\delta_{i-1>r}(r+1)\right)
     \delta(z-w).
\end{align*}

Thus

\begin{align*}
&\Big[\rho(E^1_i)(z_1),\enspace
\rho(E_i)(z_2),
\enspace \rho(E_{i-1})(w)\Big] \\
&=:a_{ii}^*a_{i-1,i}\left(\sum_{k=1}^{i-2}a_{k,i-2}a^*_{k,i-2}
     :\right) \delta(z_1-w)\delta(z_2-w)  \\
&+:a_{ii}^*a^*_{i-1,i}\left(\sum_{k=1}^{i}a_{ki}a^*_{ki}:
    -\sum_{k=1}^{i-1} a_{k,i-1}a_{k,i-1}^*\right)
     \delta(z_1-w)\delta(z_2-w) \\
&-:a_{ii}^*a^*_{i-1,i}\left(\sum_{k=1}^{i}a_{ki}a^*_{ki}:\right)
     \delta(z_1-w)\delta(z_2-w) \\
&+i\delta_{1\leq i-1\leq r}
    :a_{ii}^*(z_1)a_{ii}^*(z_2)a^*_{i-1,i-1}(z_2):\delta(z_2-w)
    \partial_{z_2}\delta(z_1-z_2) \\
&+(i+2)\delta_{1\leq i\leq r}
     a_{ii}^*(z_1)
     a^*_{i-1,i}(z_2)\delta(z_2-w)\partial_{z_2}
     \delta(z_1-z_2) \\
&+\delta_{1\leq i\leq r}
     a_{ii}^*(z_2)
     a^*_{i-1,i}(z_1)\delta(z_2-w)\partial_{z_2}
     \delta(z_1-z_2) \\
&-\delta_{1\leq i-1\leq r}
     a_{ii}^*(z_1)
     a^*_{i-1,i}(z_2)\delta(z_2-w)\partial_{z_2}
     \delta(z_1-z_2) \\
&-:a_{ii}^*
     \sum_{l=1}^{i-2}a_{l,i-2}a_{li}^*:
     \delta(z_1-w)\delta(z_2-w)  \\
&-:a_{ii}^*a_{i-1,i}^*:(b_{i-1}+b_i)\delta(z_1-z_2)
      \delta(z_2-w) \\
&-\left(\delta_{i>r}(r+1)+\delta_{i\leq r}(i+1)
     -\gamma^2\right):a_{ii}^*(z_1)a_{i-1,i}^*(z_1):
        \delta(z_2-w) \partial_w\delta(z_1-w) \\
&-\left( \gamma^2-\delta_{i-1\leq r}i-\delta_{i-1>r}(r+1)\right)
     :a_{ii}^*(z_1)
      a_{ii}^*(z_1)\partial_wa_{i-1,i-1}^*(w):
      \delta(z_1-z_2)\delta(z_2-w) \\
&+\left( \gamma^2-\delta_{i-1\leq r}i
    -\delta_{i-1>r}(r+1)\right):a_{ii}^*(z_1)
      a_{ii}^*(z_2) a_{i-1,i-1}^*(z_1):
      \partial_w\delta(z_1-w)
      \delta(z_2-w).\\
\end{align*}

Next we have

\begin{align*}
&\Big[\rho(E^2_i)(z_1), \enspace
\rho(E_i)(z_2),\enspace \rho(E_{i-1})(w)\Big] \\
&= \Bigg(-
    :a_{ii}^*(z_2)a_{i-1,i}^*(z_2)\left(
        \sum_{k=1}^{i-2}
     a_{k,i-2}(w)a_{k,i-2}^*(w)\right):  \\
&\quad\quad +:a_{ii}^*a_{i-1,i}^*
     \left(\sum_{k=1}^{i-1} a_{k,i-1}a_{k,i-1}^*\right):
     +:a^*_{i-1,i}\sum_{k=1}^{i-1}a_{k,i-1}a_{ki}^*\\
&\quad\quad+
      :a_{ii}^*\left(\sum_{l=1}^{i-2}a_{l,i-2}a_{l,i}^* \right):
    - :a_{i-1,i}^* \left(\sum_{l=1}^{i-1}a_{l,i-1}a_{li}^* \right) \\
&\qquad\quad+a_{ii}^* 
a_{i-1,i}^*b_{i-1}\Bigg)\delta(z_1-z_2)\delta(z_2-w)
\\ &-\left( \gamma^2-\delta_{i-1\leq r}i-\delta_{i-1>r}(r+1)\right)
:a_{ii}^*(z_2)a_{i-1,i}^*(z_1):
     \partial_w\delta(z_1-w)
      \delta(z_2-w).  \\
\end{align*}

The third summation contributes
\begin{align*}\Big[&\rho(E^3_i)(z_1),\enspace
\rho(E_i)(z_2),\enspace \rho(E_{i-1})(w)\Big] \\
&=:a^*_{i-1,i}a_{ii}^*b_{i}\delta(z_1-z_2)\delta(z_2-w) \\
&+\left(\delta_{i>r}(r+1)+\delta_{i\leq r}(i+1)
     -\gamma^2\right) a_{ii}^*(z_2)
     a_{i-1,i}^*(z_2)\partial_{z_1}\delta(z_1-z_2)
     \delta(z_2-w)   \\
&+\Bigg(-a_{ii}^*(z_1)a_{i-1,i}^*(z_2)
     \left( \gamma^2-\delta_{i>r}(r+1)
     -\delta_{i,r+1}\right) \\
&\quad-a_{ii}^*(z_1)a_{ii}^*(z_2)a_{i-1,i-1}^*(z_2)
         \left( \gamma^2-\delta_{i-1>r}(r+1)\right)\Bigg)
\delta(z_2-w)\partial_{z_2}\delta(z_1-z_2). \\
\end{align*}

Consequently
\begin{align*}
&\Big[\rho(E_i)(z_1),\rho(E_i)(z_2),\enspace
\rho(E_{i-1})(w)\Big] \\
&= \Big( i\delta_{1\leq i-1\leq r}
    :a_{ii}^*(z_1)a_{ii}^*(z_2)a^*_{i-1,i-1}(z_2):  \\
&-\left( \gamma^2-\delta_{i-1\leq r}i-\delta_{i-1>r}(r+1)\right)
     :a_{ii}^*(z_1)
      a_{ii}^*(z_1)\partial_wa_{i-1,i-1}^*(w): \\
&+\left( \gamma^2-\delta_{i-1\leq r}i
    -\delta_{i-1>r}(r+1)\right):a_{ii}^*(z_1)
      a_{ii}^*(z_2) a_{i-1,i-1}^*(z_1):  \\
&- \left( \gamma^2-\delta_{i-1>r}(r+1)\right)
     a_{ii}^*(z_1)a_{ii}^*(z_2)a_{i-1,i-1}^*(z_2)\Big)
    \delta(z_2-w)\partial_{w}\delta(z_1-w) \\
&+(i+2)\delta_{1\leq i\leq r}
     a_{ii}^*(z_1)
     a^*_{i-1,i}(z_2)\delta(z_2-w)\partial_{z_2}
     \delta(z_1-z_2) \\
&-\delta_{1\leq i-1\leq r}
     a_{ii}^*(z_1)
     a^*_{i-1,i}(z_2)\delta(z_2-w)\partial_{z_2}
     \delta(z_1-z_2) \\
&+\delta_{1\leq i\leq r}
     a_{ii}^*(z_2)
     a^*_{i-1,i}(z_1)\delta(z_2-w)\partial_{z_2}
     \delta(z_1-z_2) \\
&-\left(\delta_{i>r}(r+1)+\delta_{i\leq r}(i+ 1)
     -\gamma^2\right):a_{ii}^*(z_1)a_{i-1,i}^*(z_1):
        \delta(z_2-w) \partial_w\delta(z_1-w) \\
&-\left(\gamma^2-\delta_{i-1\leq r}i-\delta_{i-1>r}(r+1)\right)
:a_{ii}^*(z_2)a_{i-1,i}^*(z_1):
     \partial_w\delta(z_1-w)
      \delta(z_2-w)  \\
&-\left(\delta_{i>r}(r+1)+\delta_{i\leq r}(i+1)
     -\gamma^2\right) a_{ii}^*(z_2)
     a_{i-1,i}^*(z_2)\delta(z_2-w)\partial_{w}\delta(z_1-w) \\
&-a_{ii}^*(z_1)a_{i-1,i}^*(z_2)
     \left( \gamma^2-\delta_{i>r}(r+1)
     -\delta_{i,r+1}\right)
     \delta(z_2-w)\partial_{w}\delta(z_1-w)\\   %%%%%
&=\Bigg(\Big( (i+2)\delta_{1\leq i\leq r}
     -\delta_{1\leq i-1\leq r} -
     \gamma^2+\delta_{i>r}(r+1)
     +\delta_{i,r+1}\Big)a_{ii}^*(z_1)
     a^*_{i-1,i}(z_2)\\
&+\Big(\delta_{1\leq i\leq r}
  -\gamma^2+\delta_{i-1\leq r}i+\delta_{i-1>r}(r+1)
\Big):a_{ii}^*(z_2)a_{i-1,i}^*(z_1):
    \\
&-\left(\delta_{i>r}(r+1)+\delta_{i\leq r}(i+ 1)
     -\gamma^2\right):a_{ii}^*(z_1)a_{i-1,i}^*(z_1): \\
&-\left(\delta_{i>r}(r+1)+\delta_{i\leq r}(i+1)
     -\gamma^2\right) a_{ii}^*(z_2)
     a_{i-1,i}^*(z_2)\Bigg)
     \delta(z_2-w)\partial_{w}\delta(z_1-w)\\  \\  %%%
%&=\Big(-\left(\delta_{i>r}(r+1)+\delta_{i\leq r}(i+1)
%    -\gamma^2\right):a_{ii}^*(z_1)a_{i-1,i}^*(z_1): \\
%&-a_{ii}^*(z_1)a_{i-1,i}^*(z_2)
%    \left( \gamma^2-(i+1)\delta_{1\leq i\leq r}
%    -\delta_{i>r}(r+1) \right)\Big)
%    \delta(z_2-w)\partial_{w}\delta(z_1-w)    \\
%&-\left(\delta_{i>r}(r+1)+\delta_{i\leq r}(i+1)
%    -\gamma^2\right) a_{ii}^*(z_2) a_{i-1,i}^*(z_2)
%     \delta(z_2-w)\partial_{w}\delta(z_1-w)     \\
%&-\left(\gamma^2-\delta_{i-1\leq r}i-\delta_{i-1>r}(r+1)\right)
%:a_{ii}^*(z_2)a_{i-1,i}^*(z_1):
%    \partial_w\delta(z_1-w)
%     \delta(z_2-w)  \\
%&+\delta_{1\leq i\leq r}
%    a_{ii}^*(z_2)
%    a^*_{i-1,i}(z_1)\delta(z_2-w)\partial_{w}
%    \delta(z_1-w)                     \\  \\  %%%%%
&=\Big(-\left(\delta_{i>r}(r+1)+\delta_{i\leq r}(i+1)
     -\gamma^2\right):a_{ii}^*(z_1)a_{i-1,i}^*(z_1):  \\
&+ \left(\delta_{i>r}(r+1)+\delta_{i\leq r}(i+1)
     -\gamma^2\right)
    a_{ii}^*(z_1)a_{i-1,i}^*(z_2)  \\
&-\left(\delta_{i>r}(r+1)+\delta_{i\leq r}(i+1)
     -\gamma^2\right) a_{ii}^*(z_2) a_{i-1,i}^*(z_2)    \\
&+\left(\delta_{i>r}(r+1)+\delta_{i\leq r}(i+1)
     -\gamma^2\right) :a_{ii}^*(z_2)a_{i-1,i}^*(z_1):\Big)
      \delta(z_2-w) \partial_w\delta(z_1-w)=0. \\
\end{align*}

Now we turn to the last series of computations:

\begin{align*}
&\Big[\rho(E^1_{i-1})(z_1),\enspace
\rho(E_i)(z_2),\enspace \rho(E_{i-1})(w)\Big] \\
&=-(i-2)\delta_{i-2\leq r}
     :a_{ii}^*(z_2)a_{i-1,i-1}^*(z_1)a_{i-1,i-1}^*(z_2)
         \delta(z_2-w)\partial_{z_2}\delta(z_1-z_2)\\
%&-:a_{ii}^*a_{i-1,i-1}^* a_{i-1,i-1}^*\left(
%      \sum_{k=1}^{i-2} a_{k,i-2}a_{k,i-2}^*\right):
%      \delta(z_1-z_2)\delta(z_2-w)\\
%&+:a_{ii}^*a_{i-1,i-1}^*a_{i-1,i-1}^*\left(
%      \sum_{k=1}^{i-2} a_{k,i-2}a_{k,i-2}^*
%     \right):\delta(z_1-z_2)\delta(z_2-w) \\
&-(i+2)\delta_{i-1\leq r}:a_{ii}(z_2)a_{i-1,i-1}(z_1)a_{i-1,i-1}(z_2):
       \delta(z_2-w)\partial_{z_2}\delta(z_1-z_2)\\
&+(i-2)\delta_{i-2\leq r}
     :a_{i-1,i-1}^*(z_1)a^*_{i-1,i}(z_2):
       \delta(z_2-w)\partial_{z_2}\delta(z_1-z_2) \\
&+2\delta_{i-1\leq r}
     :a_{i-1,i-1}^*(z_1)a^*_{i-1,i}(z_2):
       \delta(z_2-w)\partial_{z_2}\delta(z_1-z_2) \\
%&+2:a_{i-1,i-1}^*a_{i-1,i-1}^*
%      \left(\sum_{l=1}^{i-1}a_{l,i-1}a_{li}^*\right)
%      \delta(z_1-z_2)\delta(z_2-w)\\
&+ 2:a_{i-1,i-1}^*
        a_{ii}^*\left(\sum_{l=1}^{i-2}a_{l,i-2}a_{l,i-1}^*
        \right)\delta(z_1-z_2)\delta(z_2-w)\\
&-:a_{i-1,i-1}^*a_{i-1,i}^*\left(
        \sum_{k=1}^{i-2} a_{k,i-2}a_{k,i-2}^*\right):
         \delta(z_1-z_2)\delta(z_2-w)\\
%&-2:a_{i-1,i-1}^*a_{i-1,i-1}^*
%      \left(\sum_{l=1}^{i-1}a_{l,i-1}a_{li}^*\right)
%      \delta(z_1-z_2)\delta(z_2-w)\\
%&- 2:a_{i-1,i-1}^*
%      a_{ii}^*\left(\sum_{l=1}^{i-2}a_{l,i-2}a_{l,i-1}^*
%      \right)\delta(z_1-z_2)\delta(z_2-w)\\
%&+2:a_{i-1,i-1}^*(z_1)a_{i-1,i}^*(z_2):
%       \delta(z_2-w)\partial_{z_2}\delta(z_1-z_2) \\
&+a_{i-1,i}^*a_{i-1,i-1}^*\left(
        \sum_{k=1}^{i-1}
      a_{k,i-1}a_{k,i-1}^*\right):
      \delta(z_1-z_2)\delta(z_2-w) \\
&-:a_{i-1,i-1}^*
        \Big(\sum_{l=1}^{i-2}a_{l,i-2}a_{li}^*:\Big)
        \delta(z_1-z_2)\delta(z_2-w)\\
&+:a_{i-1,i-1}^*a_{i-1,i-1}^*a_{ii}^*
        b_{i-1}\delta(z_1-z_2)\delta(z_2-w) \\
&-\left( \gamma^2-\delta_{i-1\leq r}i-\delta_{i-1>r}(r+1)\right)
     :a_{i-1,i-1}^*(z_1)a_{i-1,i-1}^*(z_1)a_{ii}^*(z_2):
     \delta(z_2-w)\partial_w\delta(z_1-w).
\end{align*}

Next we have

\begin{align*}
&\Big[\rho(E^2_{i-1})(z_1),\enspace
\rho(E_i)(z_2),\enspace \rho(E_{i-1})(w)\Big] \\
&=:a_{i-1,i}^*a_{i-1,i-1}^*\left(
        \sum_{k=1}^{i-2}
     a_{k,i-2}a_{k,i-2}^*\right):
     \delta(z_1-z_2)\delta(z_2-w) \\
&-2 :a_{ii}^*a_{i-1,i-1}^*
     \Big(\sum_{k=1}^{i-2}a_{k,i-2}a_{k,i-1}^*:\Big)
      \delta(z_1-z_2)\delta(z_2-w) \\
&-a_{i-1,i}^*a_{i-1,i-1}^*\left(
        \sum_{k=1}^{i-1}
      a_{k,i-1}a_{k,i-1}^*\right):
      \delta(z_1-z_2)\delta(z_2-w) \\
%&+:a^*_{i-1,i}\Big(\sum_{k=1}^{i-2}a_{k,i-2}a_{k,i-1}^*\Big)
%       \delta(z_1-z_2)\delta(z_2-w) \\
%&-:a_{i-1,i}^*\left(\sum_{l=1}^{i-2}a_{l,i-2}a_{l,i-1}^*
%   \right)
%   \delta(z_1-z_2)\delta(z_2-w)\\
&+:a_{i-1,i-1}^*\left(\sum_{l=1}^{i-2}a_{l,i-2}a_{li}^*\right)
       \delta(z_1-z_2)\delta(z_2-w)\\
&-a_{i-1,i}^*a_{i-1,i-1}^*b_{i-1}
        \delta(z_1-z_2)\delta(z_2-w)\\
&+\left( \gamma^2-\delta_{i-1\leq r}i-\delta_{i-1>r}(r+1)\right)
      a_{i-1,i}^*(z_1) \partial_wa_{i-1,i-1}^*(w)
     \delta(z_1-z_2)\delta(z_2-w). \\
\end{align*}

The third summation contributes
\begin{align*}
&\Big[\rho(E^3_{i-1})(z_1),\enspace
\rho(E_i)(z_2),\enspace \rho(E_{i-1})(w)\Big] \\
&= -:a_{ii}^*a_{i-1,i-1}^*a_{i-1,i-1}^*:b_{i-1}
     \delta(z_1-z_2)\delta(z_2-w) \\
& +\left(\delta_{i-1>r}(r+1)+\delta_{i-1\leq r}i
      -\gamma^2\right)
     :a_{ii}^*(z_2)a_{i-1,i-1}^*(z_2)a_{i-1,i-1}^*(z_2):
     \delta(z_2-w)\partial_{z_1}\delta(z_1-z_2) \\
&+:a_{i-1,i-1}^*a_{i-1,i}^*:b_{i-1}
     \delta(z_1-z_2)\delta(z_2-w) \\
&+\left(\delta_{i-1>r}(r+1)+\delta_{i-1\leq r}i
      -\gamma^2\right)
     :a_{i-1,i-1}^*(z_2)a_{i-1,i}^*(z_2):
      \delta(z_2-w)\partial_{z_1}\delta(z_1-z_2) \\
&-(\gamma^2-(r+1)\delta_{i> r+2}-\delta_{i,r+2})
     a_{i-1,i-1}^*(z_1)a_{i-1,i}^*(z_2)
     \delta(z_2-w)\partial_{z_2}\delta(z_1-z_2) \\
&+2(\gamma^2-(r+1)\delta_{i> r+1}+\frac{r}{2}\delta_{i,r+2})
     a_{i-1,i-1}^*(z_1)a_{i-1,i-1}^*(z_2)a_{i,i}^*(z_2)
     \delta(z_2-w)\partial_{z_2}\delta(z_1-z_2). \\
\end{align*}

Adding these all up we get
\begin{align*}
&\Big[\rho(E_{i-1})(z_1),\enspace
\rho(E_i)(z_2),\enspace \rho(E_{i-1})(w)\Big] \\
&=-(i-2)\delta_{i-2\leq r}
     :a_{ii}^*(z_2)a_{i-1,i-1}^*(z_1)a_{i-1,i-1}^*(z_2)
         \delta(z_2-w)\partial_{z_2}\delta(z_1-z_2)\\
%&-:a_{ii}^*a_{i-1,i-1}^* a_{i-1,i-1}^*\left(
%      \sum_{k=1}^{i-2} a_{k,i-2}a_{k,i-2}^*\right):
%      \delta(z_1-z_2)\delta(z_2-w)\\
%&+:a_{ii}^*a_{i-1,i-1}^*a_{i-1,i-1}^*\left(
%      \sum_{k=1}^{i-2} a_{k,i-2}a_{k,i-2}^*
%     \right):\delta(z_1-z_2)\delta(z_2-w) \\
&-(i+2)\delta_{i-1\leq r}:a_{ii}(z_2)a_{i-1,i-1}(z_1)a_{i-1,i-1}(z_2):
       \delta(z_2-w)\partial_{z_2}\delta(z_1-z_2)\\
&+(i-2)\delta_{i-2\leq r}
     :a_{i-1,i-1}^*(z_1)a^*_{i-1,i}(z_2):
       \delta(z_2-w)\partial_{z_2}\delta(z_1-z_2) \\
&+2\delta_{i-1\leq r}
     :a_{i-1,i-1}^*(z_1)a^*_{i-1,i}(z_2):
       \delta(z_2-w)\partial_{z_2}\delta(z_1-z_2) \\
&-\left( \gamma^2-\delta_{i-1\leq r}i-\delta_{i-1>r}(r+1)\right)
     :a_{i-1,i-1}^*(z_1)a_{i-1,i-1}^*(z_1)a_{ii}^*(z_2):
     \delta(z_2-w)\partial_w\delta(z_1-w)   \\
&+\left( \gamma^2-\delta_{i-1\leq r}i-\delta_{i-1>r}(r+1)\right)
      a_{i-1,i}^*(z_1) \partial_wa_{i-1,i-1}^*(w)
     \delta(z_1-z_2)\delta(z_2-w)  \\
& +\left(\delta_{i-1>r}(r+1)+\delta_{i-1\leq r}i
      -\gamma^2\right)
     :a_{ii}^*(z_2)a_{i-1,i-1}^*(z_2)a_{i-1,i-1}^*(z_2):
     \delta(z_2-w)\partial_{z_1}\delta(z_1-z_2) \\
&+\left(\delta_{i-1>r}(r+1)+\delta_{i-1\leq r}i
      -\gamma^2\right)
     :a_{i-1,i-1}^*(z_2)a_{i-1,i}^*(z_2):
      \delta(z_2-w)\partial_{z_1}\delta(z_1-z_2) \\
&-(\gamma^2-(r+1)\delta_{i> r+2}-\delta_{i,r+2})
     a_{i-1,i-1}^*(z_1)a_{i-1,i}^*(z_2)
     \delta(z_2-w)\partial_{z_2}\delta(z_1-z_2) \\
&+2(\gamma^2-(r+1)\delta_{i> r+1}+\frac{r}{2}\delta_{i,r+2})
     a_{i-1,i-1}^*(z_1)a_{i-1,i-1}^*(z_2)a_{i,i}^*(z_2)
     \delta(z_2-w)\partial_{z_2}\delta(z_1-z_2) \\ \\
%&=-(i-2)\delta_{i-2\leq r}
%    :a_{ii}^*(z_2)a_{i-1,i-1}^*(z_1)a_{i-1,i-1}^*(z_2)
%        \delta(z_2-w)\partial_{z_2}\delta(z_1-z_2)\\
%&-(i+2)\delta_{i-1\leq r}:a_{ii}(z_2)a_{i-1,i-1}(z_1)a_{i-1,i-1}(z_2):
%      \delta(z_2-w)\partial_{z_2}\delta(z_1-z_2)\\
%&-\left( \gamma^2-\delta_{i-1\leq r}i-\delta_{i-1>r}(r+1)\right)
%    :a_{i-1,i-1}^*(z_1)a_{i-1,i-1}^*(z_1)a_{ii}^*(z_2):
%    \delta(z_2-w)\partial_w\delta(z_1-w)   \\
%&+2(\gamma^2-(r+1)\delta_{i> r+1}+\frac{r}{2}\delta_{i,r+2})
%    a_{i-1,i-1}^*(z_1)a_{i-1,i-1}^*(z_2)a_{i,i}^*(z_2)
%    \delta(z_2-w)\partial_{z_2}\delta(z_1-z_2) \\
%& +\left(\delta_{i-1>r}(r+1)+\delta_{i-1\leq r}i
%     -\gamma^2\right)
%    :a_{ii}^*(z_2)a_{i-1,i-1}^*(z_2)a_{i-1,i-1}^*(z_2):
%    \delta(z_2-w)\partial_{z_1}\delta(z_1-z_2) \\
%&+(i-2)\delta_{i-2\leq r}
%    :a_{i-1,i-1}^*(z_1)a^*_{i-1,i}(z_2):
%      \delta(z_2-w)\partial_{z_2}\delta(z_1-z_2) \\
%&+\left(\delta_{i-1>r}(r+1)+\delta_{i-1\leq r}i
%     -\gamma^2\right)
%    :a_{i-1,i-1}^*(z_2)a_{i-1,i}^*(z_2):
%     \delta(z_2-w)\partial_{z_1}\delta(z_1-z_2) \\
%&+2\delta_{i-1\leq r}
%    :a_{i-1,i-1}^*(z_1)a^*_{i-1,i}(z_2):
%      \delta(z_2-w)\partial_{z_2}\delta(z_1-z_2) \\
%&-(\gamma^2-(r+1)\delta_{i> r+2}-\delta_{i,r+2})
%    a_{i-1,i-1}^*(z_1)a_{i-1,i}^*(z_2)
%    \delta(z_2-w)\partial_{z_2}\delta(z_1-z_2) \\
%&+\left( \gamma^2-\delta_{i-1\leq r}i-\delta_{i-1>r}(r+1)\right)
%     a_{i-1,i}^*(z_1) \partial_wa_{i-1,i-1}^*(w)
%    \delta(z_1-z_2)\delta(z_2-w)  \\  \\
&=-\left( \gamma^2-\delta_{i-1\leq r}i-\delta_{i-1>r}(r+1)\right)
     :a_{i-1,i-1}^*(z_1)a_{i-1,i-1}^*(z_1)a_{ii}^*(z_2):
     \delta(z_2-w)\partial_w\delta(z_1-w)   \\
&+2(\gamma^2-(r+1)\delta_{i> r+1}-i\delta_{i\leq r+1})
     a_{i-1,i-1}^*(z_1)a_{i-1,i-1}^*(z_2)a_{i,i}^*(z_2)
     \delta(z_2-w)\partial_{z_2}\delta(z_1-z_2) \\
& +\left(\delta_{i-1>r}(r+1)+\delta_{i-1\leq r}i
      -\gamma^2\right)
     :a_{ii}^*(z_2)a_{i-1,i-1}^*(z_2)a_{i-1,i-1}^*(z_2):
     \delta(z_2-w)\partial_{z_1}\delta(z_1-z_2) \\
&+\left(\delta_{i-1>r}(r+1)+\delta_{i-1\leq r}i
      -\gamma^2\right)
     :a_{i-1,i-1}^*(z_2)a_{i-1,i}^*(z_2):
      \delta(z_2-w)\partial_{z_1}\delta(z_1-z_2) \\
&-\left( \gamma^2-\delta_{i-1\leq r}i-\delta_{i-1>r}(r+1)\right)
     a_{i-1,i-1}^*(z_1)a_{i-1,i}^*(z_2)
     \delta(z_2-w)\partial_{z_2}\delta(z_1-z_2) \\
&+\left( \gamma^2-\delta_{i-1\leq r}i-\delta_{i-1>r}(r+1)\right)
      a_{i-1,i}^*(z_1) \partial_wa_{i-1,i-1}^*(w)
     \delta(z_1-z_2)\delta(z_2-w).  \\
     &=0 \\
\end{align*}

\end{proof}

\section{Acknowledgement}
The first author is grateful to FAPESP for financial support and to the
  University
of S\~ao Paulo for hospitality. The second author is grateful to the 
University of Sydney for
support and hospitality.

\def\cprime{$'$}
\providecommand{\bysame}{\leavevmode\hbox to3em{\hrulefill}\thinspace}
\bibliographystyle{amsplain}

\end{document}